\documentclass{amsart}
\usepackage{amssymb, amsmath, amscd, latexsym, mathrsfs, appendix}
\usepackage{graphics}
\usepackage{amsmath}
\usepackage{amsfonts}
\usepackage[all]{xy}
\newtheorem{thm}{Theorem}[section]
\newtheorem{cor}[thm]{Corollary}
\newtheorem{lem}[thm]{Lemma}

\newtheorem{prop}[thm]{Proposition}
\newtheorem{defn}[thm]{Definition}

\newtheorem{rem}[thm]{Remark}
\newtheorem{expl}[thm]{Example}

\numberwithin{equation}{section}

\newcommand{\lra}{\longrightarrow}
\newcommand{\co}{\colon\!}

\newcommand{\smin}{\smallsetminus}
\newcommand{\pt}{\star\,}  
\newcommand{\sha}{^{\sharp}}

\newcommand{\ste}{\textup{st}}

\newcommand{\id}{\textup{id}}
\newcommand{\im}{\textup{im}}

\newcommand{\holim}{\textup{holim}}
\newcommand{\hocolim}{\textup{hocolim}}

\newcommand{\map}{\textup{map}}

\newcommand{\tm}{^{\,t\!}}

\newcommand{\emb}{\textup{emb}}

\newcommand{\csigl}{\mathscr W}   

\newcommand{\lf}{{\textup{$\ell$f}}}
\newcommand{\Sm}{\textup{Sm}}
\newcommand{\intr}{\textup{int}}
\newcommand{\admi}{\mathfrak X}
\newcommand{\admib}{\admi_{\bullet}}

\newcommand{\sS}{\mathcal S}
\newcommand{\sO}{\mathcal O}

\newcommand{\NN}{\mathbb N}
\newcommand{\RR}{\mathbb R}
\newcommand{\ZZ}{\mathbb Z}

\newcommand{\twosub}[2]{\begin{array}{cc}
\scriptstyle{#1} \\  [-1mm] \scriptstyle{#2}  \end{array}}

\newcommand{\colimsub}[1]{\begin{array}[t]{cc} \textup{colim} \\
[-1.7mm] \scriptstyle{#1} \end{array}}

\newcommand{\holimsub}[1]{\begin{array}[t]{cc} \textup{holim} \\ [-1mm]
\scriptstyle{#1} \end{array}}

\newcommand{\hocolimsub}[1]{\begin{array}[t]{cc} \textup{hocolim} \\
[-1.7mm] \scriptstyle{#1} \end{array}}

\begin{document}

\title[Maps to the plane III]{Functor calculus and the discriminant method} %
\author{Rui Reis and Michael Weiss}%
\address{Mathematisches Institut, Universit\"{a}t M\"{u}nster, Einsteinstrasse 62, 48149 M\"{u}nster, Germany}%
\email{rgonc\_\,01@uni-muenster.de} \email{m.weiss\,@uni-muenster.de}  %
\thanks{This project was generously supported by the Engineering and Physical Sciences Research Council (UK),
Grant EP/E057128/1.}
\date{January 2013}%
\begin{abstract} The discriminant method is a tool for
describing the cohomology, or the homotopy type, of certain spaces of smooth maps with uncomplicated
singularities from a smooth compact manifold $L$ to $\RR^k$. We recast some of it in the
language of functor calculus. This reformulation allows us to use the discriminant method in a setting
where we wish to impose conditions on the multilocal behavior
of smooth maps $f\co L\to \RR^k$.
\end{abstract}
\maketitle
\section{Introduction} \label{sec-intro}
Vassiliev's discriminant method is a tool for
describing the cohomology, or the homotopy type, of certain spaces of smooth maps with uncomplicated
singularities from a smooth compact manifold $L$ to $\RR^k$. Here we recast some of it in the
language of functor calculus, more precisely, in the language of the manifold calculus which has been traditionally
used to study spaces of smooth embeddings \cite{WeissEmbrep,WeissEmb,GoWeEmb}. This reformulation allows
us to use the discriminant method in situations where we wish to impose conditions on the multilocal behavior
of smooth maps $f\co L\to \RR^k$, that is, the restrictions of $f$ to small neighborhoods of finite sets in $L$. \newline
We recall Vassiliev's ``main theorem'' \cite{Vassiliev1,Vassiliev2}.
Fix $r\ge 0$ and let $\admi$ be an open semi-algebraic subset of the finite dimensional vector space $P$ of polynomial maps
of degree $\le r$ from $\RR^\ell$ to $\RR^k$. This is required to be invariant under the action (by pre-composition
and truncation) of the Lie group $G$ of invertible polynomial maps $\RR^\ell\to\RR^\ell$ fixing the origin.
Next, let $L$ be a smooth manifold of dimension $\ell$, without boundary for simplicity.
A smooth map $f\co L\to \RR^k$ is considered to have an \emph{inadmissible point}, or \emph{$\admi$-inadmissible point},
at $x\in L$ if the $r$-th Taylor polynomial of $f$ at $x$, in local coordinates centered at $x$, belongs to
$P\smin \admi$. The space $C^\infty(L,\RR^k;\admi)$ of smooth maps $L\to \RR^k$ which are
everywhere admissible comes with a jet prolongation map
\begin{equation} \label{eqn-uncompjetbdl}
C^\infty(L,\RR^k;\admi) \lra \Gamma(L;\admi)~.
\end{equation}
Here $\Gamma(L;\admi)$ is the section space of a fiber bundle on $L$ with fibers homeomorphic to
$\admi$; the fiber at $x\in L$ is the space of the $r$-jets of \emph{$\admi$-admissible} map germs $(L,x)\to \RR^k$.

\begin{thm} \label{thm-vass} \cite{Vassiliev1,Vassiliev2}
If the codimension of $P\smin \admi$ in $P$ is $\ge\ell+2$ everywhere, then~{\rm (\ref{eqn-uncompjetbdl})}
induces an isomorphism in integer cohomology. If the codimension is $\ge\ell+3$ everywhere,
then~{\rm (\ref{eqn-uncompjetbdl})} is a homotopy equivalence.
\end{thm}

Vassiliev assumes that $L$ is closed. That is the difficult case; the case where
$L$ is (connected and) noncompact follows from Gromov's general $h$-principle \cite{Gromov71},\cite{Haefliger71}.
Vassiliev also has a variant where $L$ is compact with boundary, and the focus is on admissible
smooth maps $L\to \RR^k$ prescribed on/near $\partial L$. \newline
In many applications $\admi$ contains all $f\co \RR^\ell\to \RR^k$ in $P$ which are regular at the origin,
but this is not a general requirement. (For $r=0$ and $k>0$ it would not make much sense.)

\medskip
For closed $L$, let $\sO(L)$ be the poset of open subsets of $L$, ordered by inclusion.
An open set $U\in \sO(L)$ is \emph{tame} if it contains a compact codimension zero
smooth submanifold $K$ such that the inclusion of $K\smin\partial K$ in $U$ is
isotopic to a diffeomorphism. Let $\sO\tm(L)$ be the full
sub-poset consisting of the tame open sets. There are contravariant functors $\mathcal S$,
respectively $\mathcal S_{\Gamma}$,
from $\sO\tm(L)$ to chain complexes taking $U\in \sO\tm(L)$ to the singular chain complex of
$C^\infty(U,\RR^k;\admi)$, respectively of $\Gamma(U;\admi)$.
We shall apply manifold calculus to $\mathcal S$ and $\mathcal S_{\Gamma}$.
In manifold calculus it is customary to consider contravariant functors (cofunctors) from $\sO(L)$ or $\sO\tm(L)$ to the category
of \emph{spaces}. The values of $\sS$ and $\mathcal S_{\Gamma}$ are chain complexes, which we can view as spaces
by using the Kan-Dold equivalence between chain complexes graded over the positive integers and simplicial
abelian groups. Alternatively, we can take the view that manifold calculus should be applicable to
cofunctors from $\sO(L)$ or $\sO\tm(L)$ to any model category. \newline
In manifold calculus, a cofunctor $\mathcal F$ from $\sO\tm(L)$ to a model category is considered \emph{good}
if it takes isotopy equivalences in $\sO\tm(L)$ to weak equivalences \cite{WeissEmb}.
A good cofunctor on $\sO\tm(L)$ has best polynomial approximations
$\mathcal F \to T_i\mathcal F$, also known as Taylor approximations. Without unraveling the deeper meaning
of \emph{polynomial cofunctor}, which has something to do with higher excision properties,
we can define $T_i\mathcal F$ by
\[  T_i\mathcal F(U) = \holimsub{V\in \sO i(U)} \mathcal F(V) \]
where $\sO i(U)$ is the sub-poset of $\sO\tm(U)$ consisting of the open sets which are diffeomorphic to
a disjoint union of at most $i$ copies of $\RR^\ell$. In particular, $T_i\mathcal F$ is determined
up to natural weak equivalence by the restriction $\mathcal F|\sO i(L)$. We say informally that $\mathcal F$ is
\emph{analytic} if the canonical map from $\mathcal F(U)$ to $\holim_i~T_i\mathcal F(U)$
is a weak equivalence for every $U\in\sO\tm(L)$.

\begin{thm} \label{thm-vasscalc} If the codimension of $P\smin \admi$ in $P$ is $\ge\ell+2$ everywhere, then
the cofunctors $\sS$ and $\sS_{\Gamma}$ on $\sO\tm(L)$ are analytic.
\end{thm}

We now explain how theorem~\ref{thm-vass} for closed $L$ is a corollary of theorem~\ref{thm-vasscalc}.
Jet prolongation induces a natural transformation $\mathcal S \to \mathcal S_{\Gamma}$.
It is straightforward to verify that $\mathcal S(U) \to \mathcal S_{\Gamma}(U)$ is a chain equivalence if
$U$ belongs to $\sO i(L)$ for some $i$. Since both functors are analytic, this implies, by the homotopy
invariance property of homotopy inverse limits, that
$\mathcal S(U) \to \mathcal S_{\Gamma}(U)$ is a chain
equivalence for all $U\in \sO\tm(L)$. In particular, this holds for $U=L$, which proves
the first part of theorem~\ref{thm-vass}. The second part of theorem~\ref{thm-vass} follows because source and target
in~(\ref{eqn-uncompjetbdl}) are simply-connected if the codimension of $P\smin \admi$ in $P$
is $\ge \ell+3$ everywhere.

\medskip
Our proof of theorem~\ref{thm-vasscalc} uses the discriminant method and
is not far removed from Vassiliev's argument for theorem~\ref{thm-vass}. Vassiliev
constructs two spectral sequences converging to the cohomology of source and target in theorem~\ref{thm-vass},
respectively, and then proceeds by spectral sequence comparison.
We shall see that his spectral sequences are the spectral sequences associated with the Taylor towers of
$\sS(L)$ and $\sS_\Gamma(L)$.

\bigskip
Vassiliev's proof of theorem~\ref{thm-vass} suggests certain generalizations which are unfortunately
concealed in the statement. In the formulation of theorem~\ref{thm-vasscalc}
these generalizations are more accessible, as we shall now explain. \newline
On the category of finite nonempty sets and injective maps, we have a cofunctor $T\mapsto P^T$ where $P$ is
the vector space of polynomials mentioned above. Suppose given a subfunctor $\admi_\bullet$ of $T\mapsto P^T$
such that $\admi_T$ is open semi-algebraic in $P^T$, and invariant under the action of the Lie group $G^T$.
Again, let $L$ be a smooth manifold of dimension $\ell$, without boundary for simplicity.
A smooth map $f\co L\to \RR^k$ is considered to be \emph{$\admi_\bullet$-admissible}
if, for every finite nonempty subset $T$ of $L$, the $T$-tuple of $r$-th Taylor polynomials of $f$ at $t\in T$,
in any local coordinates centered at $t$, belongs to
$\admi_T$. Let $C^\infty(L,\RR^k;\admi_\bullet)$ be the space of admissible smooth maps $L\to \RR^k$
and let $\mathcal S(L)$ be its singular chain complex (recycled notation).

\begin{thm} \label{thm-vasscalcmult}
Suppose that the functor $\admi_\bullet$ is coherently open and large enough \emph{(definitions~\ref{defn-coh}
and~\ref{defn-codim})}. Then, for closed $L$, the cofunctor on $\sO\tm(L)$
defined by $U\mapsto \sS(U)$ is analytic.
\end{thm}

By general manifold calculus principles, the theorem leads to a spectral sequence converging to the homology of
$C^\infty(L,\RR^k;\admi_\bullet)$, the spectral sequence associated with the Taylor tower of $\sS(L)$.
The $p$-th column of its $E^1$-page is, broadly speaking,
the homology of the space of decorated cardinality $p$ subsets $S$ of $L$.
The decorations available for such an $S$ are in a space homeomorphic to
$\hocolim_T\admi_T$, where $T$ runs through the nonempty subsets of $S$.
For details see remark~\ref{rem-spec3}. 

\medskip
Theorem~\ref{thm-vasscalcmult} generalizes one half of theorem~\ref{thm-vasscalc}. Namely, for an
open subset $\admi$ of $P$
as in theorem~\ref{thm-vasscalc}, we define $\admib$ by $\admi_T=\admi^T\subset P^T$. Then
the functor $\sS$ of theorem~\ref{thm-vasscalcmult} coincides with the functor $\sS$ of
theorem~\ref{thm-vasscalc}. \newline

For another example which is not covered by theorem~\ref{thm-vasscalc}, but is covered by a relative version of
theorem~\ref{thm-vasscalcmult}, take
$L=D^n\times D^2$ and let the target space be the plane ($k=2$).
In~\cite{RWI}, we describe an $\admib$ which satisfies the conditions of theorem~\ref{thm-vasscalcmult}
and some additional ones which seem reasonable, among them the constraint that each $\admi_T$ can be regarded as a union
of \emph{left-right} equivalence classes, not just right-equivalence classes alias $G^T$-orbits.
When $|T|=1$, the open set $\admi_T\subset P^T$ is the
union of the left-right equivalence classes \emph{regular}, \emph{fold}, \emph{cusp}, \emph{swallowtail},
\emph{beak-to-beak} and \emph{lips}. For $|T|>1$, the open set $\admi_T$ consists of $T$-tuples
$(f_t)_{t\in T}\in P^T$ such that each $f_t$ belongs to one of the six types just mentioned, and
certain general position conditions are satisfied if, for some $y\in \RR^2$,
there are several $t\in T$ where $f_t(0)=y$ and the linear map $df_t(0)\co \RR^{n+2}\to\RR^2$
fails to be surjective. To characterize this $\admib$ as a minimal choice it was necessary in \cite{RWI}
to introduce a notion of equivalence for multigerms which is even weaker than left-right equivalence.

\medskip
We end the introduction by noting that the manifold calculus has recently been reformulated in
a way which makes closer connections with operad theory. See for example \cite{BritoWeiss} and \cite{AroneTurchin}.
In the new formulation, the functors to which the calculus is applied are typically contravariant
continuous (or \emph{enriched}) functors from the category of smooth manifolds of a fixed dimension, with smooth
embeddings as morphisms, to the category of spaces. The new framework has some conceptual advantages but the
traditional one has some practical advantages, too, where proofs are concerned. Therefore, even though our
results could be translated into the new language, for practical purposes we have opted for the traditional one.

\section{The discriminant method after Vassiliev} \label{sec-discrim}
Let $L$ be a smooth compact manifold of dimension $\ell$. We fix $k,r\ge 0$ and take $P$ to be the
vector space of polynomial maps
of degree $\le r$ from $\RR^\ell$ to $\RR^k$. As in the introduction, suppose that $\admi\subset P$ is an open
semi-algebraic subset invariant under the action of the Lie group $G$.

\begin{expl}{\rm
Popular choices for $\admi$ that satisfy the codimension condition in theorem~\ref{thm-vass}
are as follows. In the first one we have
$k=1$ and $r=3$ so that $P$ is the space of polynomial maps of degree at most $3$ from $\RR^\ell$ to $\RR$.
We let $\admi$ consist of all $f\in P$ which at the origin are either nonsingular, or have a Morse singularity,
or a birth-death singularity. This example was considered by Igusa~\cite{Igusa3,Igusa4}.\newline
In another example we have $k=2$ and $r=4$. Let $\admi$ consist of all $f\in P$ whose germ at the origin is
in one of the six left-right equivalence classes \emph{regular}, \emph{fold}, \emph{cusp}, \emph{swallowtail},
\emph{lips} and \emph{beak-to-beak} as in \cite{RWI}.
}
\end{expl}

Let $\varphi\co L\to \RR^k$ be a
smooth map such that the $r$-jets of $\varphi$ at all points
$x\in \partial L$, in any local coordinates centered at $x$, belong to $\admi$.
Let $\Sm$ be the affine
space of all smooth maps $f\co L\to \RR^k$ which satisfy
$j^rf|\partial L=j^r\varphi|\partial L$, where $j^r$ is the $r$-jet prolongation \emph{for maps
from $L$ to $\RR^k$}. This condition implies that $f$ and $\varphi$ agree on $\partial L$. If $r>0$, it is stronger since it
prescribes certain (higher) partial derivatives of $f$ involving directions normal to $\partial L$,
at points $x\in \partial L$. We equip $\Sm$ with the Whitney $C^\infty$ topology.
\newline
Let $\csigl:=C^\infty(L,\RR^k;\admi,\varphi)\subset \Sm$ be the open subset consisting of all $f$
which are $\admi$-admissible. Let $\Gamma(L;\admi,\varphi)$ be the space of sections
of the jet bundle $J^r(L,\RR^k)\to L$ which extend $j^r\varphi$ on $\partial L$ and take values in the
subbundle determined by $\admi$. We have the jet prolongation map
\begin{equation}\label{eqn-jetprolong}
\csigl=C^\infty(L,\RR^k;\admi,\varphi) \lra \Gamma(L;\admi,\varphi)~.
\end{equation}
The general (relative) form of theorem~\ref{thm-vass} is as follows.

\begin{thm} \label{thm-vassrel} \cite{Vassiliev1,Vassiliev2}
If the codimension of $P\smin \admi$ in $P$ is $\ge\ell+2$ everywhere, then~{\rm (\ref{eqn-jetprolong})}
induces an isomorphism in integer cohomology. If the codimension is $\ge\ell+3$ everywhere,
then~{\rm (\ref{eqn-jetprolong})} is a homotopy equivalence.
\end{thm}

\medskip
For the proof, we can assume that the double $L\cup_{\partial L}L$ is contained in a euclidean space $\RR^N$
as a nonsingular real algebraic subset of $\RR^N$ and $\partial L\subset L\cup_{\partial L}L$
is also a nonsingular real algebraic subset. This is justified by the Nash-Tognoli embedding
theorem. See \cite{BCR}, especially \cite[rmk~14.1.15]{BCR}.
We identify $L$ with the first summand in $L\cup_{\partial L}L$. This makes $L$ into a real semi-algebraic subset
of $\RR^N$. Indeed $L$ is the closure of the union of some connected components of the semi-algebraic set
$L\cup_{\partial L}L\smin \partial L$. Any connected component of a semi-algebraic subset of $\RR^N$
is semi-algebraic \cite[thm.2.4.5]{BCR}, and the closure of a semi-algebraic subset of $\RR^N$ is semi-algebraic
\cite[prop.2.2.2]{BCR}. \newline
By appendix~\ref{sec-boundary}, we may assume that $\varphi\co L\to \RR^k$ extends to a polynomial map $\RR^N\to \RR^k$.
Also by appendix~\ref{sec-boundary}, we may assume that $r$ (in the description of $P$ and $\admi\subset P$) is
even and strictly positive.

\begin{thm}\label{thm-interpol}
There exists an ascending sequence
$(A_i)_{i\in \NN}$ of finite dimensional (dimension $d_i$) affine
subspaces of $\Sm$ with the following properties.
\begin{itemize}
\item[(i)] \emph{Dense}: the union of the $A_i$ is dense in $Sm$;
\item[(ii)] \emph{Algebraic}: every $f\in A_i$ extends to a polynomial map on $\RR^N$;
\item[(iii)] \emph{Tame}: for $f\in  A_i$~, the number of $x\in L$ where $f$ is $\admi$-inadmissible
is bounded above by a constant $\alpha_i\in \NN$~;
\item[(iv)] \emph{Interpolating}: for every $i$ and every $T\subset L\smin\partial L$ with $|T|\le i$,
the projection $A_i\to \prod_{x\in T} J^r_x(\RR^\ell,\RR^k)\cong P^T$ is onto.
\end{itemize}
\end{thm}

For a proof of this, see Appendix~\ref{sec-interpol}. This follows Vassiliev in
all essentials, but we have eliminated some complicated transversality arguments.

\medskip
By the density property, the inclusion
\[ \colimsub{i}(A_i\cap \csigl) \lra \csigl  \]
is a weak homotopy equivalence. Hence we can approximate the cohomology of $\csigl$ with the
cohomology of $A_i\cap \csigl $. The reduced cohomology $\tilde H^s(A_i\cap \csigl)$ is isomorphic to the
locally finite homology
\[  H_{d_i-s-1}^\lf(A_i\smin \csigl) \]
by a form of Poincar\'{e}-Alexander-Lefshetz duality in the finite dimensional affine space $A_i$.
See appendix~\ref{sec-exci}.
Here we are using the fact that
\[  B_i= A_i\smin \csigl \]
is an ENR, euclidean neighborhood retract, an essential condition
for the duality statement. Indeed, $B_i$ is a closed semialgebraic set in $A_i$. (Triangulation
theorems in \cite{BCR} for example imply that closed semialgebraic sets in euclidean spaces are ENRs.)
Still following Vassiliev, we therefore focus on $B_i$ and its locally finite homology. (This explains
the expression \emph{discriminant method}, since $B_i\subset A_i$ can be called the discriminant variety.)

\medskip
We construct a ``resolution''
\begin{equation} \label{eqn-reso} RB_i\to B_i~,
\end{equation}
a proper homotopy equivalence of locally compact spaces, where $RB_i$ admits a filtration
which will help us to understand its homological properties. In detail, $RB_i$ is the classifying space of
a topological poset whose elements are pairs $(f,T)$ where $f\in B_i$ and $T\subset L$ is a
\emph{bad event}, i.e., a finite nonempty set of inadmissible points for $f$. (We note
that such a $T$ satisfies $T\subset L\smin\partial L$ since $B_i\subset A_i\subset\Sm$.)
The order relation is
given by $(f,S)\le (g,T)$
iff $f=g$ and $S\subset T$. There is an obvious metrizable topology on the underlying set
where we say that a sequence $((f_j,S_j))_{j\in\NN}$ converges to $(f,S)$ if $(f_j)_{j\in\NN}$ converges to $f$,
in the Whitney topology on $C^\infty(L,\RR^k)$,
and $(S_j)_{j\in\NN}$ converges to $S$ in the Hausdorff topology. See appendix~\ref{sec-badtop}.
With this topology, the order
relation ``$\le$'' is closed. This leads to a canonical
choice of topology on the classifying
space of the topological poset (see appendix~\ref{sec-badtop} again). \newline
The fiber of the projection
\[  RB_i \lra B_i \]
at $f\in B_i$ is the classifying space of the finite poset of all subsets of $L\smin\partial L$
which are bad events (hence finite and nonempty) for $f$.
Since that poset has a maximal element, we learn from this that~(\ref{eqn-reso})
has contractible and compact fibers.
By appendix~\ref{sec-badtop}, it is also a proper map of locally compact spaces.

\begin{lem} \label{lem-reso} The space $RB_i$ is an ENR and the
resolution map $RB_i\to B_i$ is a proper homotopy equivalence.
\end{lem}

\proof It is shown in appendix~\ref{sec-badtop} that $RB_i$ is an ENR. The fibers of $RB_i\to B_i$ are
simplices. This implies that $RB_i\to B_i$ is a \emph{cell-like} map according to \cite{Lacher69},
and so by \cite[Thm 1.2]{Lacher69} it is a proper homotopy equivalence.
\qed

\medskip
We filter $RB_i$ as follows. The underlying topological poset is filtered such that $(f,S)$ lives in the
$s$-th stage of the poset, where $s=|S|$. This determines a filtration
\[ F_1 RB_i \subset F_2 RB_i \subset F_3 RB_i \subset \cdots \subset RB_i \]
such that all points in the interior of a nondegenerate simplex determined by a diagram
\[  (f,S_0) < (f,S_1) < \cdots < (f,S_{n-1}) < (f,S_n) \]
belong to $F_s RB_i$, where $s=|S_n|$.

\begin{lem} \label{lem-resres} Each $F_pRB_i$ is an ENR.
The restricted resolution map $F_iRB_i \to B_i$ induces an isomorphism in locally finite homology
in dimensions greater than $d_i-i$.
\end{lem}

\proof It is shown in appendix~\ref{sec-badtop} that $F_pRB_i$ is an ENR. Now fix $i>0$.
Let $C\subset B_i$ consist of all $f\in B_i$ which admit a bad event
of cardinality $\ge i$. Let $T=\{1,2,\dots,i\}$ and let
\[  Z\subset A_i\times\emb(T,L\smin\partial L) \]
consist of all $(f,e)$ such that $e(T)$ is a bad event for $f$.
This has codimension $\ge i(\ell+2)$ in $A_i\times\emb(T,L\smin\partial L)$ by property (iv)
in theorem~\ref{thm-interpol} (see~\ref{rem-codimestimate}), and
so has codimension $\ge 2i$ relative to $A_i$. The set $Z$ is a semialgebraic set
and it follows that its image in $A_i$, which is $C$, is also semialgebraic; see \cite[Prop 2.2.7]{BCR}. \newline
Now let $F_iRB_i|C$ be the portion of $F_iRB_i$ projecting to $C$.
Then $F_iRB_i|C$ is the image of a map
\[  Z \times \Delta^{i-1} \lra F_iRB_i \]
where $\Delta^{i-1}$ plays the role of classifying space of the poset $\{1,2,3,\dots,i\}$.
(By appendix~\ref{sec-badtop}, this map can also be interpreted as a semi-algebraic map.)
The codimension of $Z \times \Delta^{i-1}$ relative to $A_i$ is $\ge i+1$. Therefore
the codimension of $F_iRB_i|C$ relative to $A_i$ is also $\ge i+1$. As the projection
$F_iRB_i|C\to C$ is onto by definition, it follows that the
codimension of $C$ in $A_i$ is $\ge i+1$. \newline
Next, there is a commutative square
\begin{equation} \label{eqn-MVloc}
\xymatrix@M5pt{ F_iRB_i|C \ar[r] \ar[d] & F_iRB_i \ar[d] \\
C \ar[r] & B_i
}
\end{equation}
By appendix~\ref{sec-badtop}, it is a square of ENRs and proper maps. Let $B_i'$ be the
pushout of
\[
\xymatrix@M5pt{ C & \ar[l] F_iRB_i|C \ar[r]  & F_iRB_i
}
\]
so that there is a canonical map $B'_i\to B_i$. As $F_iRB_i|C\to F_iRB_i$ is a cofibration, it follows that
$B'_i$ is an ENR; see \cite{Hu}.
The fibers of $B'_i\to B_i$ over points of $C$ are singletons, and the other fibers are
standard simplices. This implies that $B'_i\to B_i$ is a cell-like map,
and so by \cite{Lacher69} it is a proper homotopy equivalence. Therefore it induces an isomorphism
\begin{equation} \label{eqn-locfinres} H_*^\lf(B'_i)\cong H_*^\lf(B_i)~.
\end{equation}
>From the construction of $B'_i$ and the comparison with $B_i$~, the square~(\ref{eqn-MVloc}) determines a long exact
Mayer-Vietoris sequence relating the groups
\[  H_*^\lf(B_i),\qquad  H_*^\lf(F_iRB_i)\oplus H_*^\lf(C),\qquad H_*^\lf(F_iRB_i|C).  \]
With our (co-)dimension estimates this proves that
$H_*^\lf(F_iRB_i) \to  H_*^\lf(B_i)$
is an isomorphism when $*>d_i-i$. \qed

\medskip
Summing up the main insights so far, we have
\[ \tilde H^s(A_i\cap \csigl)~\cong~H_{d_i-s-1}^\lf(B_i)~\cong~H_{d_i-s-1}^\lf(F_iRB_i) \]
if $s+1<i$. Therefore we shall examine the locally finite
homology of $F_iRB_i$ and use the filtration
\[ F_1 RB_i \subset F_2 RB_i \subset F_3 RB_i \subset \cdots \subset F_iRB_i \]
for that. By appendix~\ref{sec-exci}, the locally finite homology of the pair $(F_p RB_i,F_{p-1}RB_i)$ is
isomorphic to the locally finite homology of
\[  F_p RB_i\smin F_{p-1}RB_i \]
since $F_{p-1}RB_i$ is a closed sub-ENR of $F_pRB_i$.
Elements of $F_p RB_i\smin F_{p -1}RB_i$ have the form $(f,S,x)$ where $f\in A_i$ and $S\subset L\smin\partial L$
is a bad event for $f$ with $|S|=p$, while $x$ is an element of the simplex spanned by $S$.
(The coordinate $x$ is in the interior of the simplex.) There is a map
\begin{equation}\label{eqn-vassilspec}
 F_p RB_i\smin F_{p -1}RB_i \lra \binom{L\smin\partial L}{p}
 \end{equation}
whose target is the space of unordered configurations of $p$ points in $L\smin\partial L$.
It is defined by $(f,S,x)\mapsto S$.

\begin{thm}\label{thm-loctriv}
The map~(\ref{eqn-vassilspec}) is a locally trivial projection, for $0<p\leq i$. The fibers
are homeomorphic to $\RR^b\times\RR^{p-1}\times (P\smin\admi)^p$ where $b=d_i-p\dim(P)$.
\end{thm}

\proof
To show that it is locally trivial we factorize the map as follows:
\begin{equation} \label{eqn-vassilspecfac}
F_p RB_i\smin F_{p-1}RB_i \lra E_p  \lra  \binom{L\smin\partial L}{p}.
\end{equation}
Here $E_p$ is the space of triples $(h,S,x)$ where $S$ is an unordered configuration of $p$ points
in $L\smin\partial L$,
\[ h\in \prod_{s\in S} J^r_s(L,\RR^k) \]
has coordinates $h_s$ which are all inadmissible (i.e., not in $\admi$ when local coordinates centered at $s$
are used) and $x$ is an element in the interior of the simplex spanned by $S$. The left-hand arrow associates
to $(f,S,x)$ as above the triple $(h,S,x)$, where $h=j_S^rf$. The right-hand arrow is forgetful, and it is
clearly a fiber bundle projection with fibers homeomorphic to $(P\smin\admi)^p\times \RR^{p-1}$.
It remains to show that the left-hand arrow in~(\ref{eqn-vassilspecfac})
is also locally trivial. Given $(h,S,x)\in E_p$~, the portion of $F_p RB_i\smin F_{p-1}RB_i$ being
mapped to $(h,S,x)$ is the affine subspace
\begin{equation} \label{eqn-affinewonder} \{f\in A_i~|~j^r_Sf=h~\}~.
\end{equation}
Most important is the observation that the affine space~(\ref{eqn-affinewonder}) is nonempty
and has a dimension which is independent of $(h,S,x)$. Indeed, it is a translate of the kernel
of the projection
\[   A_i\lra \prod_{s\in S} J^r_s(L,\RR^k)~\cong~P^S \]
which is assumed to be onto by the interpolation condition. The condition is applicable
because $|S|=p \le i$.
\qed

\bigskip
We now ask how the spectral sequence
in locally finite homology determined by the filtrations
\[  F_1 RB_i \subset F_2 RB_i \subset F_3 RB_i \subset \cdots \subset F_iRB_i \]
depends on $i$. (Note that theorem~\ref{thm-loctriv} describes the $E^1$ page of this
spectral sequence and shows that the $E^1$ page is independent of $i$ in the ``active''
region.) The affine quotient $A_{i+1}/A_i$ is a vector space, assuming $A_i\ne \emptyset$.
For $p\le i$ there is a map
\begin{equation} \label{eqn-suspend} \Psi\co F_p RB_{i+1} \lra A_{i+1}/A_i
\end{equation}
defined by $(f,S,x)\mapsto [f]$ in the notation of the proof of theorem~\ref{thm-loctriv}. We do not claim that $\Psi$
is a bundle projection, but we still wish to have the homological corollaries.

\begin{prop} \label{prop-speccompat}
The spectral sequences in locally finite homology determined by the filtrations
\[  F_1 RB_i \subset F_2 RB_i \subset F_3 RB_i \subset \cdots \subset F_iRB_i \]
and
\[  F_1 RB_{i+1} \subset F_2 RB_{i+1} \subset F_3 RB_{i+1} \subset \cdots \subset F_iRB_{i+1} \]
are isomorphic up to a shift of $d_{i+1}-d_i$, equal to the dimension of $A_{i+1}/A_i$.
\end{prop}

\proof First we look at the fibers
of~(\ref{eqn-suspend}). The fiber $\Psi^{-1}(0)$ is exactly $F_p RB_i$. The fiber $\Psi^{-1}(v)$ is like $F_p RB_i$
but constructed from the translate $A_i+\bar v$ of $A_i$ in $Sm$, where $\bar v\in A_{i+1}$ represents $v$. Now $A_i+\bar v$
satisfies roughly the same conditions that $A_i$ satisfies (e.g.~the interpolation condition), so that we have
a result for $\Psi^{-1}(v)$ identical to~\ref{thm-loctriv} for $\Psi^{-1}(0)=F_p RB_i$. \newline
More generally, let us choose a triangulation of $A_{i+1}/A_i$ by linear simplices. Then the above arguments show that,
for simplices $\sigma$ and $\tau$ in $A_{i+1}/A_i$ with $\sigma\subset\tau$, the inclusion
\[   \Psi^{-1}(\sigma) \lra \Psi^{-1}(\tau) \]
induces an isomorphism in locally finite homology. It follows immediately that there is a Leray-Serre type spectral sequence
converging to the locally finite homology of the source $F_p RB_{i+1}$ of $\Psi$, with second page equal to the locally
finite homology of $A_{i+1}/A_i$ with coefficients in the locally finite homology of the fiber of $\Psi$. For obvious reasons
this collapses and so we obtain
\[   H_*^\lf(F_p RB_{i+1}) \cong H_{*-d_{i+1}+d_i}^\lf(F_p RB_i). \]
A similar relative argument shows
\[   H_*^\lf(F_p RB_{i+1},F_s RB_{i+1}) \cong H_{*-d_{i+1}+d_i}^\lf(F_p RB_i,F_s RB_i) \]
for $s<p$. Therefore by exact couple technology, the two spectral sequences in the lemma are isomorphic
up to a shift as stated. \qed

\begin{rem} \label{rem-belatedspec}
{\rm The above proposition~\ref{prop-speccompat} provides, among other things, a belated justification for the
strategy which we adopted earlier in this section by claiming or assuming that $H^*(A_i\cap \csigl)$ is a
good approximation to $H^*(\csigl)$. Now we can be more precise. There is a commutative ladder of homomorphisms
\begin{equation}\label{eqn-csiglladder}
\xymatrix@R=12pt{
\vdots \ar[d]   &  \vdots  \\
\tilde H^*(A_{i+1}\cap\csigl)  \ar[d]   & H_{d_{i+1}-*-1}^\lf(F_{i+1}RB_{i+1}) \ar[l] \ar[u] \\
\tilde H^*(A_{i}\cap\csigl)  \ar[d]   & H_{d_{i}-*-1}^\lf(F_iRB_{i}) \ar[l] \ar[u] \\
\tilde H^*(A_{i-1}\cap\csigl)  \ar[d]   & H_{d_{i-1}-*-1}^\lf(F_{i-1}RB_{i-1}) \ar[l] \ar[u] \\
\vdots & \vdots \ar[u]
}
\end{equation}
The maps in the left-hand column are induced by the inclusions $A_i\to A_{i+1}$, and the maps
in the right-hand column are as in proposition~\ref{prop-speccompat}. The horizontal map in row $i$ is the
composition
\[
\xymatrix{
H_{d_{i}-*-1}^\lf(F_iRB_{i}) \ar[r] & H_{d_{i}-*-1}^\lf(B_{i}) \ar[r]^-{\cong} & \tilde H^*(A_{i}\cap\csigl)
}
\]
By lemma~\ref{lem-resres}, it is an isomorphism for $*+1<i$. By the dimension formula in theorem~\ref{thm-loctriv},
the arrow
\[
H_{d_{i-1}-*-1}^\lf(F_{i-1}RB_{i-1}) \lra H_{d_{i}-*-1}^\lf(F_iRB_{i})
\]
in the right-hand column of the ladder is also an isomorphism when $*+1<i$. (Any deviation must be reflected in the
locally finite homology of
$F_iRB_i\smin F_{i-1}RB_i$. By theorem~\ref{thm-loctriv}, this is the total space of a bundle and as such
has dimension
\[
\le~i\ell+(d_i-i\dim(P))+i-1+i(\dim(P)-\ell-2)~=~d_i-i-1
\]
where the summand $i\ell$ is the dimension of the base space of the bundle.)
It follows that the arrow
\[  \tilde H^*(A_i\cap\csigl) \lra  \tilde H^*(A_{i-1}\cap\csigl) \]
is an isomorphism for $*<i-2$. \newline
It follows also
that the spectral sequence described in proposition~\ref{prop-speccompat}
(stabilized with respect to $i$) \emph{converges} to the reduced cohomology of $\csigl$.
Setting this up as a homological spectral sequence in the fourth quadrant, we have
\[   E^1_{p,q}=H_{p+q+d_i-1}^\lf(F_pRB_i\smin F_{p-1}RB_i)  \]
for some or any $i\ge p$. The differentials have the form
\[   d^s\co E^s_{p,q} \lra E^s_{p-s,q+s-1}~. \]
The dimension of $F_pRB_i\smin F_{p-1}RB_i$ is $\le d_i-p-1$,
by theorem~\ref{thm-loctriv} again.
Therefore all $E^1_{p,q}$ where $p+q+d_i-1>d_i-p-1$ are zero. This means
$E^1_{p,q}= 0 \textup{ for } q>-2p$
which leads to the following picture of the $E^1$ page (fat dots for
potentially nonzero positions):
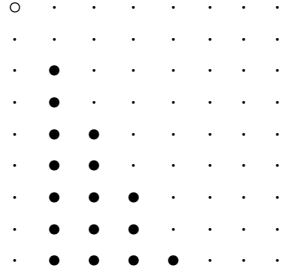
\[
\begin{matrix}
\circ &  \cdot & \cdot & \cdot & \cdot & \cdot & \cdot & \cdot & \\
\cdot &  \cdot & \cdot & \cdot & \cdot & \cdot & \cdot & \cdot & \\
\cdot &  \bullet & \cdot & \cdot & \cdot & \cdot & \cdot & \cdot & \\
\cdot &  \bullet & \cdot & \cdot & \cdot & \cdot & \cdot & \cdot & \\
\cdot &  \bullet & \bullet & \cdot & \cdot & \cdot & \cdot & \cdot & \\
\cdot &  \bullet & \bullet & \cdot & \cdot & \cdot & \cdot & \cdot & \\
\cdot &  \bullet & \bullet & \bullet & \cdot & \cdot & \cdot & \cdot & \\
\cdot &  \bullet & \bullet & \bullet & \cdot & \cdot & \cdot & \cdot & \\
\cdot &  \bullet & \bullet & \bullet & \bullet & \cdot & \cdot & \cdot &

\end{matrix}
\]
By reflecting at the origin, marked by a circle $\circ$~, we can also view this as a
cohomology spectral sequence in the second quadrant.
}
\end{rem}

\bigskip
This completes our (re)construction of Vassiliev's spectral sequence converging to $H^*(\csigl)$. Vassiliev uses the same
ideas to construct an analogous spectral sequence converging to the reduced cohomology of
$\Gamma(L;\admi,\varphi)$. Jet prolongation induces a map between the two spectral sequences which
specializes to an isomorphism of the $E^1$ pages. This implies that the jet prolongation map itself
induces an isomorphism from $H^*(\Gamma(L;\admi,\varphi))$ to $H^*(\csigl)$. \newline
We shall not follow this path. Instead we will finish this section with a few observations on finiteness
properties of $H^*(\csigl)$ and the Vassiliev spectral sequence. Then in section~\ref{sec-discman} we shall use
the spectral sequence directly to prove the first half of theorem~\ref{thm-vasscalc}.
The proof of the second half uses an independent argument which is easier. Theorem~\ref{thm-vass} can then be
deduced from theorem~\ref{thm-vasscalc} as explained in the introduction.

\begin{rem}[Finite generation] \label{rem-fin} {\rm Let $X$ be a closed semi-algebraic set in some $\RR^n$. Then $H_*^\lf(X)$
is finitely generated as a graded abelian group. To show this we identify the one-point compactification $X^\omega$ of
$X$ with the closure in $\RR^{n+1}$ of $f(X)\subset \RR^{n+1}$, where $f\co \RR^n\to S^n\subset \RR^{n+1}$ is the inverse of
stereographic projection:
\[  f(x)= sx+(0,0,\dots,0,1-s) \in S^n~,\qquad s=\frac{2}{1+\|x\|^2}~. \]
As $f$ is an algebraic map, $f(X)$ is a semi-algebraic set and so is its closure in $\RR^{n+1}$, by \cite{BCR}.
Therefore $X^\omega$ is a compact ENR and so we have
\[  H_*^\lf(X) \cong H_*^\ste(X^\omega,\omega) \cong H_*(X^\omega,\omega)~ \]
where $H_*^\ste$ is Steenrod homology; see appendix~\ref{sec-exci}. The isomorphism between Steenrod homology and ordinary
homology in this case comes from the fact that Steenrod homology satisfies the seven Eilenberg-Steenrod axioms, which
characterize ordinary homology on compact ENRs. \newline
The consequence for us is that $B_i$ in the above description of Vassiliev's spectral sequence has
finitely generated $H_*^\lf(B_i)$. Also, $RB_i$ and the subspaces $F_pRB_i$ were constructed as geometric
realizations of simplicial semi-algebraic sets, with proper face and degeneracy operators, and satisfying
an upper bound condition due to the tameness property in theorem~\ref{thm-interpol}. It follows that
$H_*^\lf(RB_i)$ and $H_*^\lf(F_pRB_i)$ are finitely generated as graded abelian groups. Hence the
spectral sequence described in proposition~\ref{prop-speccompat} and remark~\ref{rem-belatedspec} has an $E^1$-page
which is finitely generated in each bi-degree. Consequently $H^*(\csigl)$ is finitely generated in each degree.
}
\end{rem}

\section{Discriminant method and manifold calculus} \label{sec-discman}
Let $\sO(L)$ be the poset of open subsets of $L$ which contain $\partial L$, ordered by inclusion.
An open set $U\in \sO(L)$ is \emph{tame} if it contains a compact codimension zero
smooth submanifold $K$, containing $\partial L$, such that the inclusion of $\breve{K}$ in $U$ is
isotopic to a diffeomorphism (relative to $\partial L$), where $\breve{K}=(K\smin\partial K)\cup \partial L$. We write
\[  \sO\tm(L)=\{U\in\sO(L)~|~U\textup{ is tame }\}. \]
Fix $r\ge 0$ and $\admi\subset P$ as in section~\ref{sec-discrim}. The definitions
of $C^\infty(L,\RR^k;\admi,\varphi)$ and $\Gamma(L;\admi,\varphi)$
can be extended in a straightforward manner to make
contravariant functors on $\sO(L)$ by
\begin{eqnarray}
\label{eqn-S} \sS(U) & = & \textup{singular chain complex of }C^\infty(U,\RR^k;\admi,\varphi) \\
\label{eqn-SGamma} \sS_\Gamma(U) & = & \textup{singular chain complex of }\Gamma(L;\admi,\varphi).
\end{eqnarray}
The general form of theorem~\ref{thm-vasscalc} is as follows.

\begin{thm} \label{thm-vasscalcrel} If the codimension of $\admi$ is $\ge \ell+2$, then
the cofunctors $\sS$ and $\sS_\Gamma$ on $\sO\tm(L)$ defined by~(\ref{eqn-S}) and~(\ref{eqn-SGamma}) are analytic.
\end{thm}

The definition of \emph{analytic} used here generalizes that given in the introduction.
For $i\ge 0$ there is a subposet $\sO i(L)\subset\sO\tm(L)$. It consists of the $U\in \sO\tm(L)$
which can be written as a disjoint union of an open collar on $\partial L$ and an open subset of $L\smin\partial L$
abstractly diffeomorphic to a disjoint union of $\le i$ copies of $\RR^\ell$. By saying that
a good cofunctor $F$ on $\sO\tm(L)$ is analytic we mean that the canonical map
\[ F(U) \lra  \holimsub{i}\holimsub{U\supset V\in\sO i} F(V) \]
is a weak homotopy equivalence.

\bigskip
For $U$ in $\sO\tm(L)$, let $\kappa(U)$ be the poset, ordered by inclusion, of compact codimension zero
smooth submanifolds of $U$ containing $\partial L$. Let
\[  \csigl(U)=C^\infty(U,\RR^k;\admi,\varphi) \]
so that $\csigl(L)=\csigl$. Also for $K\in \kappa(L)$ let
\[  \csigl_K=\{ f\in \Sm~|~f\textup{ admissible on }K\}~. \]

\begin{lem} \label{lem-compactlim} For $U\in \sO\tm(L)$, there is a chain of natural (in $U$) homotopy equivalences
\[ \csigl(U)\simeq \cdots \simeq \holimsub{K\in\kappa(U)} \csigl_K~. \]
There is also a chain of natural homotopy equivalences between the singular cochain complex of
$\csigl(U)$ and the homotopy direct limit, over $K\in \kappa(U)$, of the singular
cochain complexes of the $\csigl_K$.
\end{lem}

\proof For every $K\in \kappa(U)$, the restriction map $\csigl_K\to C^\infty(K,\RR^k;\admi,\varphi)$
is a homotopy equivalence. (This uses a smooth form of the Tietze-Urysohn extension principle which
goes back to Borel; see \cite[4.9]{Broecker}.) Therefore it is enough to show that the map
\begin{equation} \label{eqn-tameres}
C^\infty(U,\RR^k;\admi,\varphi) \simeq \holimsub{K\in\kappa(U)} C^\infty(K,\RR^k;\admi,\varphi)
\end{equation}
is a homotopy equivalence.
Now let $K\in \kappa(U)$ be such that the inclusion $\breve{K}\to U$ is isotopic in $U$ to a diffeomorphism
relative to $\partial L$. Let
\[ (h_t\co \breve{K} \to U)_{t\in[0,1]} \]
be such an isotopy, so that $h_0\co \breve{K} \to U$ is the inclusion.
Then we have a map
\[  C^\infty(K,\RR^k;\admi,\varphi) \lra C^\infty(U,\RR^k;\admi,\varphi) \]
given by composing with the map $U\to K$ obtained by inverting $h_1$. This is easily seen to be a homotopy
inverse for the restriction map
\[ C^\infty(U,\RR^k;\admi,\varphi) \lra C^\infty(K,\RR^k;\admi,\varphi). \]
As the set $\kappa'(U)$ of those $K\in \kappa(U)$ for which the inclusion $\breve{K}\to U$ is isotopic in $U$
to a diffeomorphism
(relative to $\partial L$) is cofinal in $\kappa(U)$, it follows that~(\ref{eqn-tameres}) is a homotopy
equivalence. \newline
The same reasoning shows that we have a chain of natural homotopy equivalences relating the
singular chain complex of $C^\infty(U,\RR^k;\admi,\varphi)$ to the homotopy inverse limit,
over $K\in\kappa(U)$, of the singular chain complexes of $\csigl_K$. It does not matter much whether
we take that homotopy inverse limit over $\kappa(U)$ or over $\kappa'(U)$, as $\kappa'(U)$ is cofinal
in $\kappa(U)$. Since, for $K_1\le K_2\in \kappa'(U)$, the restriction map
\[ \csigl_{K_2}\to \csigl_{K_1} \]
is a homotopy equivalence, the induced map of singular chain or cochain complexes is also a homotopy
equivalence. It follows that
\[  \begin{array}{rcl}
&& \hom\Big(\holimsub{K\in \kappa(U)}\big(\textup{singular chain complex of }\csigl_K\big),\ZZ\Big) \\
&\simeq & \hocolimsub{K\in \kappa(U)}\big(\textup{singular cochain complex of }\csigl_K\big)~. \qquad\qquad \qed
\end{array}
\]

\medskip
In the following we apply the discriminant method to the spaces $\csigl_K$ where possible, and draw
conclusions about $\csigl(U)$ by means of lemma~\ref{lem-compactlim}.
We still have the affine spaces $A_i\subset \Sm$ and
the ENR subspaces $B_i=A_i\smin\csigl$ from section~\ref{sec-discrim}. For $K\in\kappa(L)$, let
\[  B_{i,K}=A_i\smin \csigl_K~. \]
Then $B_{i,K}$ is a closed subset of $B_i$. It consists of all $f\in A_i$ which have inadmissible points in $K$.
There is probably no reason to think that it is always an ENR, but in the proof of the following lemma we identify
many cases when it is.

\begin{lem}\label{lem-compactAdual} There is a chain of natural homotopy equivalences relating
the homotopy colimit, over $K\in \kappa(U)$, of the cochain complexes
\[ C_{d_i-*-1}^\lf(B_{i,K}) \]
to the homotopy colimit, over $K\in\kappa(U)$, of $\tilde C^*(\csigl_K\cap A_i)$.
\end{lem}

\proof To clarify notation, the degree 5 part of $C_{d_i-*-1}^\lf(B_{i,K})$
is $C_{d_i-6}^\lf(B_{i,K})$, for example. We speak of a \emph{cochain} complex because differentials raise
degree by one. \newline
It is not claimed that we have
\[ C_{d_i-*-1}^\lf(B_{i,K})~\simeq~ \tilde C^*(\csigl_K\cap A_i) \]
for \emph{every} $K\in\kappa(U)$. Suppose however that $K$
happens to be a semi-algebraic subset of $L$. Then $B_{i,K}$ is a semi-algebraic subset of $A_i$ and
therefore an ENR. By appendix~\ref{sec-exci} there is a chain of natural homotopy equivalences
\[  C_{d_i-*-1}^\lf(B_{i,K}) \lra \cdots \longleftarrow \tilde C^*(\csigl_K\cap A_i)~. \]
It remains to note that the set of those $K$ in $\kappa(U)$ which are semi-algebraic
is cofinal in $\kappa(U)$. This is easy. For arbitrary $K\in\kappa(U)$, choose a smooth function
$f\co K\to \RR$ having $0$ as a regular value with preimage $\partial K\smin\partial L$, and such that
$K$ is the preimage of $[0,\infty)$ under $f$. This $f$ can be
approximated by a polynomial function $g\co K\to\RR$. If the approximation is sharp enough (on values and
first derivatives), then the preimage of $[0,\infty)$ under $g$
is an element of $\kappa(U)$ close to $K$ in the sense that a small isotopy inside $U$ (rel $\partial L$)
will move it to $K$. \qed

\medskip
Let us now fix $K$ in $\kappa(L)$. Recall that there are a resolution $RB_i\to B_i$ and a filtration of $RB_i$ by
subspaces $F_pRB_i$. Let $RB_{i,K}\to B_{i,K}$ be defined similarly: $RB_{i,K}$ is the classifying space
of a topological poset whose elements are pairs $(f,S)$ where $f\in A_i$ and $S$ is a bad event
for $f$ \emph{contained in $K$}. Clearly $RB_{i,K}$ is a subspace of $RB_i$; beware that it is not defined as the preimage
of $B_{i,K}$ under the resolution map $RB_i\to B_i$. Also, we can construct a filtration $F_pRB_{i,K}$ of $RB_{i,K}$
by $F_pRB_{i,K}:=F_pRB_{i}\cap RB_{i,K}$. There is a pullback square
\begin{equation}\label{eqn-pullbackconf}
\xymatrix{ F_pRB_{i,K}\smin F_{p-1}RB_{i,K} \ar[r] \ar[d] & F_pRB_{i}\smin F_{p-1}RB_{i} \ar[d] \\
{\binom{K\smin \partial L}{p}} \ar[r] & {\binom{L\smin \partial L}{p}}}
\end{equation}
where the horizontal maps are inclusions and the right-hand vertical map is the one from theorem~\ref{thm-loctriv}.

\begin{prop}[Omnibus] \label{prop-omnibus1}
Let $U\in \sO\tm (L)$. There is a cofinal subposet $\kappa^{alg}(U)$ of $\kappa(U)$ such that the following
are satisfied for each $K\in \kappa^{alg}(U)$ and every $i>0$:
\begin{itemize}
\item[(a)] The resolution map $RB_{i,K}\to B_{i,K}$ is a proper homotopy equivalence of ENRs.
\item[(b)] Each $F_pRB_{i,K}$ is an ENR.
\item[(c)] For $p\le i$, the restricted resolution map
$F_pRB_{i,K}\to B_{i,K}$ induces an isomorphism in locally
finite homology in degrees greater than $d_i-p$.
\item[(d)] The locally finite homology of $B_{i,K}$, $RB_{i,K}$ and $F_pRB_{i,K}$ is finitely generated
as a graded abelian group.
\end{itemize}
\end{prop}

\proof The cofinal subposet $\kappa^{alg}(U)$ is implicit in the proof of lemma~\ref{lem-compactAdual}. It consists
of all $K\in \kappa(U)$ such that $K$ is a semi-algebraic set in $L$. The proof of (a) is like the proof
of lemma~\ref{lem-reso}: $RB_{i,K}$ and $B_{i,K}$ are both ENRs, the fibers of $RB_{i,K}\to B_{i,K}$ are simplices, and so
$RB_{i,K}\to B_{i,K}$ is a proper homotopy equivalence by \cite{Lacher69}. The proofs of (b) and (c) are like the proof
of lemma~\ref{lem-resres}. The proof of (d) is as in remark~\ref{rem-fin}.
\qed

\begin{rem} \label{rem-Kladder} {\rm
Given $K\in\kappa^{alg}(U)$, there is a commutative ladder of homomorphisms
\begin{equation}\label{eqn-csiglKladder}
\xymatrix@R=12pt{
\vdots \ar[d]   &  \vdots  \\
\tilde H^*(A_{i+1}\cap\csigl_K)  \ar[d]   & H_{d_{i+1}-*-1}^\lf(F_{i+1}RB_{i+1,K}) \ar[l] \ar[u] \\
\tilde H^*(A_{i}\cap\csigl_K)  \ar[d]   & H_{d_{i}-*-1}^\lf(F_iRB_{i,K}) \ar[l] \ar[u] \\
\tilde H^*(A_{i-1}\cap\csigl_K)  \ar[d]   & H_{d_{i-1}-*-1}^\lf(F_{i-1}RB_{i-1,K}) \ar[l] \ar[u] \\
\vdots & \vdots \ar[u]
}
\end{equation}
It is analogous to (\ref{eqn-csiglladder}).
By property (c) of the Omnibus proposition, the horizontal map in row $i$ is an isomorphism for $*+1<i$.
By the pullback square (\ref{eqn-pullbackconf}),
the arrow
\[
H_{d_{i-1}-*-1}^\lf(F_{i-1}RB_{i-1,K}) \lra H_{d_{i}-*-1}^\lf(F_iRB_{i,K})
\]
in the right-hand column of the ladder is also an isomorphism when $*+1<i$. It follows that the arrow
\[  H^*(A_i\cap\csigl_K) \lra  H^*(A_{i-1}\cap\csigl_K) \]
is an isomorphism for $*<i-2$.
}
\end{rem}

\bigskip
For an integer $p\ge 1$, let $E_{i,p}$ be the covariant functor from $\sO\tm(L)$ to cochain
complexes taking $U\in \sO\tm(L)$ to
\[  \colimsub{K\in\kappa(U)} C_{d_i-*-1}^\lf(F_pRB_{i,K})~. \]
We also write $E_{i,\infty}$ for the (monotone) union or direct limit of the $E_{i,p}$~,
and $E_{i,0}=0$ by convention.
Since all cochain maps used in constructing the above colimit are cofibrations, i.e., degreewise split
injective, the colimit can also be regarded as a homotopy colimit. There are cofibrations
$E_{i,p}\subset E_{i,p+1}$ for $p\ge 0$.

\begin{thm} \label{thm-Efilthomogeneous} For $p$ with $1\le p\le i$, the covariant
functor $E_{i,p}/E_{i,p-1}$ is homogeneous of degree $p$.
\end{thm}

\proof We know already that $E_{i,p}/E_{i,p-1}$ can be identified (up to natural chain equivalence)
with the covariant functor on $\sO\tm(L)$ given by
\begin{equation} \label{eqn-Efilt}
U \mapsto \colimsub{K\in\kappa(U)} C_{d_i-*-1}^\lf(F_pRB_{i,K}\smin F_{p-1}RB_{i,K})~.
\end{equation}
By~(\ref{eqn-vassilspec}) and theorem~\ref{thm-loctriv} and the pullback square~(\ref{eqn-pullbackconf}),
there is a bundle projection
\begin{equation} \label{eqn-compactloctriv}
F_pRB_{i,K}\smin F_{p-1}RB_{i,K} \lra \binom{K\smin\partial L}{p}
\end{equation}
whose fibers are described in theorem~\ref{thm-loctriv}. Taking colimits over $K\in \kappa(U)$,
we still have a bundle projection
\begin{equation} \label{eqn-betterEfilt}
 F_pRB_{i,U}\smin F_{p-1}RB_{i,U} \lra \binom{U\smin\partial L}{p}
\end{equation}
where
\[ F_pRB_{i,U}=\bigcup_{K\in\kappa(U)}F_pRB_{i,K}~. \]
Therefore $E_{i,p}/E_{i,p-1}$
can be described as the functor taking $U$ to the (re-indexed) complex of locally finite singular chains in the
total space of~(\ref{eqn-betterEfilt}) \emph{supported over a compact set of the symmetric product
$U^p/\Sigma_p$}.
To show that this functor is homogeneous of degree $p$, we need to show that
\begin{itemize}
\item[(i)] it is polynomial of degree $\le p$, according to a definition of \emph{polynomial} for covariant
functors which is given in appendix~\ref{sec-polydual};
\item[(ii)] it takes every object in $\sO m(L)$, with $m<p$, to a weakly contractible chain complex.
\end{itemize}
The ``true reason'' why these properties hold is that the functor~(\ref{eqn-Efilt}), in the description
just given, has the
(dual of the) standard form of a homogeneous degree $p$ functor given for example
in \cite{WeissEmb}. But rather than trying to make a conversion, we give a direct argument. \newline
Suppose therefore that $U$ comes equipped with
pairwise disjoint closed subsets $C_0,\dots,C_p$ in $U\smin\partial L$, as in appendix~\ref{sec-polydual}.
We assume that the $C_t$ are pairwise disjoint tame co-handles. Let $S=\{0,1,\dots,p\}$.
For $T\subset S$ let $C_T=\bigcup_{t\in T}C_t$.
In order to establish (i) we need to show that the $S$-cube of chain complexes
\[  T\mapsto \frac{E_{i,p}}{E_{i,p-1}}(U\smin C_T) \]
is cocartesian. This is equivalent to saying that the canonical map
\[ \hocolimsub{\emptyset\ne T\subset S} \frac{E_{i,p}}{E_{i,p-1}}(U\smin C_T) \lra
\frac{E_{i,p}}{E_{i,p-1}}(U) \]
is a weak equivalence. Here we may adopt the interpretation suggested above, so that
the target is the complex of locally finite singular chains in the
total space of~(\ref{eqn-betterEfilt}) supported over a compact set of
$U^p/\Sigma_p$. Then each
$(E_{i,p}/E_{i,p-1})(U\smin C_T)$ can be identified with a chain subcomplex of the target, consisting
of the locally finite singular chains supported over a compact set of $(U\smin C_T)^p/\Sigma_p$. The
homotopy colimit can be replaced by the internal sum of chain subcomplexes, because
the system of these chain subcomplexes is Reedy cofibrant~\cite{Hirschhorn}.
Therefore we have to show that the
inclusion
\begin{equation}\label{eqn-Efilt2} \sum_{\emptyset\ne T\subset S} \frac{E_{i,p}}{E_{i,p-1}}(U\smin C_T) \lra
\frac{E_{i,p}}{E_{i,p-1}}(U)
\end{equation}
is a weak homotopy equivalence. The target is still the complex of locally finite singular chains in the
total space of~(\ref{eqn-betterEfilt}) with support over a compact set of $U^p/\Sigma_p$. The source
is the subcomplex of the target generated by all locally finite singular chains which,
for some $t$, have support in a compact set of $(U\smin C_t)^p/\Sigma_p$. \newline
That~(\ref{eqn-Efilt2}) is a weak equivalence follows immediately from a subdivision argument \cite[III.7]{Dold}
and the identity
\[  \binom{U\smin \partial L}{p} = \bigcup_{\emptyset\ne T\subset S}\binom{U\smin(C_T\cup\partial L)}{p}. \]
This proves (i). \newline
For property (ii), we suppose that $U\in \sO m(L)$ for some $m<p$. This means that $U$ is the disjoint union of
an open collar on $\partial L$ and at most $m$ (without loss of generality, exactly $m$)
open subsets of $L$ diffeomorphic to $\RR^\ell$. As we have to
approximate $U$ by compact subsets, we may focus on a compact $K$ which is the disjoint union of a compact
collar and $m$ compact subsets of $L$ diffeomorphic to $D^\ell$. Then it suffices
to show that the locally finite homology of
\[ F_pRB_{i,K}\smin F_{p-1}RB_{i,K} \]
is zero. For that it is enough, by~\ref{prop-flasque},
to show that the locally compact space $F_pRB_{i,K}\smin F_{p-1}RB_{i,K}$ can be written in the form
$Y\times(0,1]$ where $Y$ is another locally compact space. For that it is enough, by local triviality
of~(\ref{eqn-compactloctriv}),
to show that every connected component of
\[ \binom{K\smin\partial K}{p} \]
can be written in that form. We write $K\cong \partial L\times[0,1]~\cup~D^\ell\times\{1,2,\dots,m\}$.
By selecting a connected component of the configuration space we are selecting a function
\[  f\co \{0,1,2,\dots,m\} \lra \NN=\{0,1,2,\dots\} \]
which counts how many points of the configuration are in each component of $K$. In particular $f(0)$ counts
the number of points from the configuration in the collar $\partial L\times[0,1]$. Because $m<p$,
we have either $f(0)>0$ or $f(s)>1$ for some $s$ in $\{1,\dots,m\}$. These cases need to be treated
separately. If $f(0)>0$, we make a function
\[ \binom{K\smin\partial K}{p} \lra (0,1] \]
by taking a configuration $S$ to the maximum of all $t\in (0,1]$ for which there is an element of $S$
which belongs to the collar $\partial L\times[0,1]$ and has second coordinate $t$. It is easy to show
that the function is the projection of a product onto a factor. If $f(s)>1$ for some $s\in \{1,\dots,m\}$,
we make a function
\[ \binom{K\smin\partial K}{p} \lra (0,1] \]
by taking a configuration $S$ to the maximum of all $t\in(0,1]$ for which there is an element of $S$
which belongs to $D^\ell\times s$ and whose distance (in standard euclidean coordinates for the disk)
from the center of the disk is $t$. Again it is easy to show that the function is the
projection of a product onto a factor.
\qed

\bigskip
\proof[Proof of first part of~\ref{thm-vasscalcrel}]
It suffices to show that $\sS(U)\to T_p\sS(U)$ is approximately
$p$-connected for every $U\in\sO\tm(L)$. We fix $p>0$ and some $i\gg p$. \newline
We ought to look at $\sS(U)$, the singular chain complex of $\csigl(U)$ as a contravariant
functor of $U$. Because of lemmas~\ref{lem-compactlim} and~\ref{lem-compactAdual} it is easier for us to
work for a little while with the reduced singular cochain complex of $\csigl(U)$ as a \emph{covariant}
functor of $U$. We have the following diagram of natural transformations between covariant functors
\[
\xymatrix@R=18pt{
{\tilde C^*(\csigl(U))} & &   {E_{i,p}(U)} \ar[d] \\
{\rule{0mm}{5mm}\hocolimsub{K} \tilde C^*(\csigl_K)} \ar[u]_-{\simeq} \ar[r]
& {\rule{0mm}{5mm}\hocolimsub{K}\tilde C^*(\csigl_K\cap A_i)}
\ar[r]^-{\simeq} & {E_{i,\infty}(U)}
}
\]
where $K$ runs through $\kappa(U)$.
We abbreviate this to
\[
\xymatrix{ Y & \ar[l]_-{\simeq} {Y^\sharp}  \ar[r] & {Y_i} \ar[r]^-{\simeq} & {E_{i,\infty}} & \ar[l] {E_{i,p}}
}
\]
suppressing the variable $U$. We apply the duality functor $D$ (see appendix~\ref{sec-polydual}) to obtain
\[
\xymatrix{ DY \ar[r]^-{\simeq} & {DY^\sharp} & \ar[l]{DY_i} & \ar[l]_-{\simeq}  {DE_{i,\infty}} \ar[r] &{DE_{i,p}}
}
\]
and form the Taylor approximations $T_p$~,
\[
\xymatrix@R=15pt{ DY \ar[r]^-{\simeq} \ar[d] & {DY^\sharp} \ar[d]
& \ar[l]{DY_i} \ar[d] & \ar[l]_-{\simeq}  {DE_{i,\infty}} \ar[r] \ar[d] &{DE_{i,p}} \ar[d]^-{\simeq} \\
T_pDY \ar[r]^-{\simeq} & {T_pDY^\sharp} & \ar[l]{T_pDY_i} &
\ar[l]_-{\simeq}  {T_pDE_{i,\infty}} \ar[r] &{T_pDE_{i,p}}~.
}
\]
In the top row, all maps are approximately $i$-connected by remarks~\ref{rem-Kladder} and~\ref{rem-belatedspec},
except the one on the right which is approximately $p$-connected by the same remarks.
In the bottom row, all maps except the one on the right are therefore approximately $(i-\ell p-p)$-connected by lemma~\ref{lem-polydim}
in appendix~\ref{sec-polydual}. Noting that $T_pDE_{i,p}\simeq DE_{i,p}\simeq T_pDE_{i,i}$,
we can view the one on the right as being induced by the approximately $i$-connected map
$DE_{i,\infty}\to DE_{i,i}$. It is therefore also approximately $(i-\ell p-p)$-connected by lemma~\ref{lem-polydim}.
It follows that the left-hand
vertical arrow is approximately $p$-connected. Since $DY$ is the reduced form of $\sS$, this completes the proof. \qed

\bigskip\bigskip
It remains to prove the second part of theorem~\ref{thm-vasscalcrel}, which states that $\sS_\Gamma$ on $\sO\tm(L)$
is analytic. In appendix~\ref{sec-polydual}, we have collected some definitions and facts related to the
notions of analytic and polynomial functor. We proceed through a few examples where these will be used.

\begin{expl} \label{expl-mapfunctor} {\rm Let $Y$ be a space, homotopy equivalent to a CW-space. Let $\psi$ be
any map from $\partial L$
to $Y$. We define $F$ on $\sO\tm(L)$ by
\[  F(U)=\map(U,Y;\psi), \]
the space of maps  $U\to Y$ which extend $\psi$. It is well known that $F$ is polynomial of degree $\le 1$, in the
sense that
\[  F(U) \lra \holimsub{W\in \sO 1(U)} F(W) \]
is a homotopy equivalence. Curiously this does not immediately tell us for which $\rho$ and $c$ the functor
$F$ is $\rho$-analytic with excess $c$ (definition~\ref{defn-analyticrc}).
It \emph{is} clear
that, in the situation of definition~\ref{defn-analyticrc} and with our choice of functor $F$,
the cube $F(U\smin C_\bullet)$ is cartesian alias $\infty$-cartesian
if $j>0$. So only the case $j=0$ is of interest. In that case we write $C$ for $C_0$ (the unique co-handle)
and $q$ for its codimension. The cube consists of a single map
\[  F(U) \lra F(U\smin C) \]
and we hope to be able to specify $\rho$, $c$ in such a way that it is always $(c+\rho-q)$-connected.
Now suppose for simplicity that $Y$ is $\ell$-connected. Then we can show that $F$ is $(\ell+1)$-analytic
with excess $0$. Indeed, there is a homotopy pullback square
\[
\xymatrix{
F(U) \ar[r] \ar[d] & F(U\smin C) \ar[d] \\
\map(D^q,Y) \ar[r] &  \map(S^{q-1},Y)
}
\]
because, up to homotopy equivalence, $U$ is obtained from $U\smin C$ by attaching a $q$-cell. The
lower horizontal map is $(\ell+1-q)$-connected, and so the upper horizontal map
is also $(\ell+1-q)$-connected. This is what we had to show. \newline
This calculation is not useless. It follows from~\cite[2.1]{WeissHomEmb} that if a good cofunctor $F$ from
$\sO\tm(L)$ to spaces is $\rho$-analytic with excess $c\ge 0$, then $\lambda F$ is also
$\rho$-analytic with excess $c$, where $\lambda F(U)$ is the singular chain complex of $F(U)$.
We emphasize that $\lambda F$ is typically not polynomial of degree $\le 1$
because the singular chain complex functor fails to preserve finite homotopy inverse limits
(such as products and homotopy pullbacks). It is typically not polynomial of any degree.
In the situation above where $F(U)=\map(U,Y;\psi)$ and $Y$ is $\ell$-connected, the Taylor tower of $\lambda F(U)$
converges to $\lambda F(U)$ for every $U\in \sO\tm(L)$, since every $U\in\sO\tm(L)$ admits a tame handle
decomposition (relative to a collar on $\partial L$) with handles of index $\le \ell$ only.
For $U=L$, we obtain therefore a second quadrant spectral
sequence converging to the ``homotopy groups'' of $\lambda F(L)$ which in this case we may
interpret as the homology groups of $F(L)$. Some details on the $E^1$-page of this spectral sequence
can be found in remark~\ref{rem-spectralE1}.
}
\end{expl}

\begin{expl} \label{expl-secfunctor} {\rm There is a variant of example~\ref{expl-mapfunctor} where we start
with a fibration $E\to L$ and a section $\psi$ of $E|\partial L$. Let $F$ be the cofunctor on $\sO\tm(L)$
defined by
\[  F(U)=\Gamma(E|U;\psi) \]
(space of sections of $E|U$ which extend $\psi$ on $\partial L$). If the fibers of $E\to L$ are $\ell$-connected,
then $F$ is $(\ell+1)$-analytic with excess $0$ and so $\lambda F$ is also $(\ell+1)$-analytic with excess $0$.
It follows that the Taylor tower of $\lambda F(U)$
converges to $\lambda F(U)$ for every $U\in \sO\tm(L)$. For $U=L$, this leads to a second quadrant spectral
sequence converging to the homology groups of $F(L)$.
}
\end{expl}

\proof[Proof of second part of theorem~\ref{thm-vasscalcrel}] The functor $\sS_\Gamma$ has the form $\lambda F$ described in
example~\ref{expl-secfunctor}, for a fiber bundle $E\to L$ whose fibers are homeomorphic
to $\admi$. The codimension of $P\smin\admi$ in the vector space $P$ is $\ge \ell+2$ everywhere, and so
$\admi$ is $\ell$-connected. Hence $\sS_\Gamma=\lambda F$ is $(\ell+1)$-analytic with excess $0$.
The Taylor tower of $\sS_\Gamma(U)$ therefore converges to $\sS_\Gamma(U)$ for every $U\in \sO\tm(L)$. \qed

\medskip
\begin{rem} \label{rem-spectralE1} {\rm The spectral sequence described in example~\ref{expl-mapfunctor},
converging to the graded group $\pi_*(\lambda F(L))\cong H_*(\map(L,Y;\psi))$, has the form
\begin{equation} \label{eqn-spectralE1}
\begin{array}{ccc} E^1_{-p,q} & = & H_c^{-q}\big(\binom{L\smin\partial L}{p}; \mathbf U \big)
\end{array}
\end{equation}
where $H_c^*$ is generalized cohomology with compact supports and
$\mathbf U$ is a fibered spectrum. The fiber of $\mathbf U$ over a configuration $S\subset L\smin  \partial L$,
where $|S|= p$, is the spectrum
\[ \mathbf H\ZZ \wedge (Y^{*S}*S^0)~. \]
Here $\mathbf H\ZZ$ is the Eilenberg-MacLane spectrum and
$Y^{*S}\cong Y*Y*\cdots *Y$ is the $S$-fold join power of $Y$. We use the base point of
$S^0$ as a base point for the join $Y^{*S}*S^0$. The space $Y^{*S}*S^0$ can also be described
as the total homotopy cofiber of the $S$-cube of spaces
$T\mapsto Y^T$, where $T\subset S$. \newline
If $Y$ is $\ell$-connected, which we are assuming, then $Y^{*S}*S^0$ is $(\ell p+2p-1)$-connected and so that cohomology group is trivial
if $q\le(\ell p+2p-1)-\ell p$, that is, if $q<2p$.
This leads to the following picture of the $E^1$ page (fat dots for
potentially nonzero positions):
\[
\begin{matrix}
\cdot & \cdot &  \cdot & \cdot & \bullet & \bullet & \bullet & \cdot   \\
\cdot & \cdot &  \cdot & \cdot & \bullet & \bullet & \bullet & \cdot   \\
\cdot & \cdot &  \cdot & \cdot & \cdot & \bullet & \bullet & \cdot   \\
\cdot & \cdot &  \cdot & \cdot & \cdot & \bullet & \bullet & \cdot   \\
\cdot & \cdot &  \cdot & \cdot & \cdot & \cdot & \bullet & \cdot   \\
\cdot & \cdot &  \cdot & \cdot & \cdot & \cdot & \bullet & \cdot   \\
\cdot & \cdot &  \cdot & \cdot & \cdot & \cdot & \cdot & \cdot   \\
\cdot & \cdot &  \cdot & \cdot & \cdot & \cdot & \cdot & \bullet   \\
\end{matrix}
\]
Formula~(\ref{eqn-spectralE1}) is a consequence of the description in \cite{WeissEmb} of the
$p$-th homogeneous layer of a good cofunctor with values in based spaces or spectra, in our case
the cofunctor $\lambda F$. To use this
we need to view $\lambda F$ as a cofunctor with spectrum values, using the stable Kan-Dold
construction to transform chain complexes into spectra. See appendix~\ref{sec-polydual}, proof of lemma~\ref{lem-polydim}.
Then, according to \cite{WeissEmb}, the $p$-th homogeneous layer of $\lambda F$, evaluated on $L$, is
the spectrum of sections with compact support of a certain fibered spectrum $\mathbf V$ on the configuration space
\[  \binom{L\smin\partial L}{p}~. \]
The fiber of $\mathbf V$ over a configuration $S$ is the total homotopy fiber of the cube
\[  T\mapsto \lambda F(N(T)) \]
where $T\subset S$ and $N(T)$ is a tubular neighborhood of $T\cup\partial L$ in $L\smin\partial L$. Here $\lambda F(N(T))$ should be
regarded as a spectrum and this simplifies to $\mathbf H\ZZ\wedge(Y^T)_+$. In the category of spectra,
total homotopy fibers of $S$-cubes agree with total homotopy cofibers up to a shift by $|S|$, so that
$\mathbf V=\Omega^p\mathbf U$. Consequently $\pi_{q-p}$ of the $p$-th homogeneous layer of $\lambda F(L)$ is
identified with
\[ \begin{array}{ccc}
H_c^{p-q}\big(\binom{L\smin\partial L}{p}; \mathbf V \big) & \cong &
H_c^{-q}\big(\binom{L\smin\partial L}{p}; \mathbf U \big).
\end{array}
\]
}
\end{rem}

\section{Admissible multijets} \label{sec-multdiscrim}
We now make the assumptions described just before theorem~\ref{thm-vasscalcmult}. In particular, for each finite nonempty set
$T$ we have an open semi-algebraic subset $\admi_T\subset P^T$, invariant under $G^T$. The assignment $T\mapsto
\admi_T$ defines a subfunctor $\admib$ of $T\mapsto P^T$.

\begin{defn} \label{defn-badness}
{\rm Let $f=(f_t)_{t\in T}\in P^T$. A nonempty subset $S\subset T$ is a
\emph{minimal bad event} for $f$ if $(f_t)_{t\in S}\notin \admi_S$ and
for every proper nonempty subset $R\subset S$ we have $(f_t)_{t\in R}\in \admi_R$.
A nonempty subset $S\subset T$ is a
\emph{bad event} for $f=(f_t)_{t\in T}$ if it is a union of minimal bad events.
The \emph{complexity} of such a bad event $S$ is the
maximum of the integers $j$ such that there exists a chain of bad subevents
\[  S_0\subsetneq S_1\subsetneq \dots \subsetneq S_{j-1}\subsetneq S_j=S~. \]
Evidently the complexity of a bad event $S$ for $f=(f_t)_{t\in T}$ is not greater than $|S|$.
}
\end{defn}

\begin{expl} \label{expl-cheap} {\rm A cheap way to construct $\admib$ as above is to specify an
open semi-algebraic
subset $\admi\subset P$, as in theorems~\ref{thm-vass} and~\ref{thm-vasscalc}, and to define
\[  \admi_T:= \admi^T\subset P^T~. \]
In this situation, for $f\in \Sm$, a subset $T$ of $L\smin\partial L$ is a bad event if and only if
each point $t\in T$ is inadmissible for $f$. A bad event $T$ is minimal if and only if $|T|=1$.\newline
For another example, where the minimal bad events can have cardinality as large as 3, see  \cite[ch.~4]{RWI}.
}
\end{expl}

\begin{defn} \label{defn-coh} {\rm The functor $\admib$ is \emph{coherently open} if the following holds.
Suppose given an open $U\subset\RR^\ell$ and a finite nonempty $S\subset U$ and a smooth $f\co U\to \RR^k$ such that
the multijet $j^r_Sf$ belongs to $\admi_S$. Then there exist a neighborhood $V\subset U$ of $S$ and a
neighborhood $W$ of $f$ in $C^\infty(U,\RR^k)$
such that for every $g\in W$ and finite nonempty subset $T$ of $V$, the multijet $j^r_Tg$ belongs to $\admi_T$. }
\end{defn}

\begin{rem} {\rm Without loss of generality, $V$ in definition~\ref{defn-coh} is a tubular neighborhood of $S$.
Then every finite nonempty $T\subset V$ comes with a preferred map $\lambda_T$ to $S$ which takes $t\in T$
to the unique $s\in S$ in the same component of $V$. For those $T$ where $\lambda_T$ is injective, the coherence
condition does not add anything to what we already have by assuming that $\admib$ is an open subfunctor of $T\mapsto P^T$.
The interesting cases are those where $\lambda_T$ is not injective. \newline
The functor $\admib$ constructed in \cite[ch.~4]{RWI}, where $\RR^k=\RR^2$, is coherently open. This is not
explicitly verified in \cite{RWI}, but the proof is easy due to another property of $\admib$ which holds by
construction: each $\admi_S\subset P^S$ is a union of finitely many left-right equivalence classes.
Consequently, for prescribed $|S|$, only a finite list of germs $f\co (U,S)\to \RR^2$ with $j^r_Sf\in \admi_S$
has to be tested for the condition in definition~\ref{defn-coh}.
}
\end{rem}

\begin{defn} \label{defn-codim} {\rm The functor $\admib$ is \emph{large enough} if the following holds.
\begin{itemize}
\item For every finite nonempty set $T$,
the elements of $P^T$ having a bad event of cardinality $p$ constitute a subset of codimension at
least $p\ell+2$ everywhere.
\item There is an upper bound $\zeta$ for the cardinality of minimal bad events.
\item There is a function $(p,j)\mapsto c(p,j)\in\NN$, defined on pairs of integers $(p,j)$ where
$p\ge 1$ and $0\le j<p$, such that
\[  \lim_{p\to\infty} c(p,j)-p\ell-j= \infty \]
and for every finite nonempty set $T$,
the elements of $P^T$ having a bad event of cardinality $p$ and
complexity $j$ constitute a subset of codimension at least $c(p,j)$ everywhere. 
\end{itemize}
}
\end{defn}

\begin{rem} {\rm Condition (i) in definition~\ref{defn-codim} ensures that the space
of all smooth $\admib$-admissible maps from $L$ to $\RR^k$ (without any boundary conditions for the moment)
is nonempty and even path connected. Namely, multijet transversality theorems imply that the set of
$f\in C^\infty(L,\RR^k)$
whose $s$-fold $r$-th order multijet prolongation (a section of a bundle on the
configuration space of $s$-element subsets of $L$) makes a transverse intersection with
the subbundle associated to
\[ P^s\smin\admi_{\{1,\dots,s\}}~, \]
for every $s\ge 0$, is dense in $C^\infty(L,\RR^k)$.
But a transverse intersection is an empty intersection by the codimension
condition in (i); therefore such an $f$ is admissible.
More generally, the set of continuous paths $(f_t\co L\to \RR^k)_{[t\in[0,1]}$ in $C^\infty(L,\RR^k)$
whose $s$-fold $r$-th order multijet prolongation makes a transverse intersection with
the subbundle associated to $P^s\smin\admi_{\{1,\dots,s\}}$~,
for every $s\ge 0$, is dense in the space of such paths. This remains
true if we only allow paths $(f_t)$ with fixed admissible $f_0$ and $f_1$.
Here again, a transverse intersection is an empty intersection by the codimension
condition in (i). Therefore each $f_t$ in a path satisfying the transversality condition is admissible. \newline
Condition (ii) is there to ensure convergence of a spectral
sequence generalizing Vassiliev's spectral sequence described in section~\ref{sec-discrim}.
For any positive integer $p$, we define $\theta(p)=\min_{j<p}~(c(p,j)-p\ell-j)$. Then the graph of
\begin{equation} \label{eqn-vanish}
 p~\mapsto~2-p-\theta(p)
\end{equation}
is a vanishing curve in the $E^1$-page of the spectral sequence.
In the situation of section~\ref{sec-discrim} alias example~\ref{expl-cheap}, where
$\admi_T=\admi^T\subset P^T$ and $\admi$ has codimension $\ge \ell+2$ in $P$, we can take $c(p,j)=p(\ell+2)$.
Then the function~(\ref{eqn-vanish}) turns into $p\mapsto 1-2p$. Its graph is the vanishing line which
we already saw in remark~\ref{rem-belatedspec}.
}
\end{rem}

The following technical definition is of interest only when $\partial L$ is nonempty.

\begin{defn} \label{defn-pure} {\rm Let $T$ be a finite nonempty set, $(f_t)_{t\in T}\in P^T$
and $S\subset T$ a bad event for $f$. The coordinates $f_t$ with $t\in S$ are elements of $P$.
We say that they \emph{participate} in the bad event $S$ for $(f_t)_{t\in T}$. \newline
An element $g$ of $P$ is considered \emph{pure} or \emph{$\admib$-pure} if it has a neighborhood
in $P$ consisting of elements which do not participate in any bad event.
}
\end{defn}

We now assume that a smooth map $\varphi\co L\to \RR^k$ has been selected which is
$\admib$-pure at every point in $\partial L$ (that is, for every $x\in \partial L$, the
jet $j^r_x\varphi$ is $\admib$-pure, in local coordinates about $x$).
This condition on $\varphi$ is rather strong, for general $\admib$ and
$L$ with nonempty boundary.

\begin{expl} {\rm For $\admib$ as in the first part of example~\ref{expl-cheap}, a map $\varphi\co L\to \RR^k$
is $\admib$-pure at all points of $\partial L$ if and only if it is $\admib$-admissible
at all points of $\partial L$.
For $\admib$ as in \cite[ch.~4]{RWI}, a map $\varphi\co L\to \RR^k$
is $\admib$-pure at all points of $\partial L$ if and only if it is regular at all
points of $\partial L$. This means that $D_x\varphi\co T_xL\to \RR^k$ is onto
for every $x\in \partial L$. It does not mean that $D_x\co T_x\partial L\to \RR^k$
is onto for every $x\in \partial L$.
}
\end{expl}

By analogy with section~\ref{sec-discrim}, let $\Sm$ be the affine
space of all smooth maps $f$ from $L$ to $\RR^k$ which satisfy
$j^rf|\partial L=j^r\varphi|\partial L$, where $j^r$ is the $r$-jet prolongation for maps
from $L$ to $\RR^k$. We equip $\Sm$ with the Whitney $C^\infty$ topology.

\begin{defn} \label{defn-mapbadness} {\rm
For $f\in \Sm$ and a nonempty finite subset $T$ of $L$, we say that $T$ is a (minimal) \emph{bad event}
for $f$ if $T$ is a (minimal) bad event for the multijet $j^r_Tf$, expressed in local
coordinates centered at the points of $T$. (This implies $T\subset L\smin\partial L$.)
}
\end{defn}

Let $\csigl:=C^\infty(L,\RR^k;\admib,\varphi)\subset \Sm$ be the open subset consisting of all $f$
which are $\admib$-admissible. Similarly, for every $U\in \sO\tm(L)$ we have
$\csigl(U):=C^\infty(U,\RR^k;\admib,\varphi)$. Let $\sS(U)$ be the singular chain complex of $\csigl(U)$.
The general (relative) form of theorem~\ref{thm-vasscalcmult} is as follows.

\begin{thm} \label{thm-vasscalcmultrel} If $\admib$ is large enough and coherently open, then
the functor $\sS$ is analytic.
\end{thm}

For the proof, we assume as in sections~\ref{sec-discrim} and~\ref{sec-discman}
that the double $L\cup_{\partial L}L$ is contained in a euclidean space $\RR^N$
as a nonsingular real algebraic subset of $\RR^N$ and $\partial L\subset L\cup_{\partial L}L$
is also a nonsingular real algebraic subset.
We identify $L$ with the first summand in $L\cup_{\partial L}L$.
We may assume that $\varphi\co L\to \RR^k$ extends to a polynomial map $\RR^N\to \RR^k$,
and that $r$ is even and strictly positive.

\begin{thm}\label{thm-interpol2}
There exists an ascending sequence
$(A_i)_{i\in \NN}$ of finite dimensional (dimension $d_i$) affine
subspaces of $\Sm$ with the following properties.
\begin{itemize}
\item[(i)] \emph{Dense}: the union of the $A_i$ is dense in $Sm$;
\item[(ii)] \emph{Algebraic}: every $f\in A_i$ extends to a polynomial map on $\RR^N$;
\item[(iii)] \emph{Tame}: for $f\in  A_i$~, the cardinality of any bad event for $f$
is bounded above by a constant $\alpha_i\in \NN$~;
\item[(iv)] \emph{Interpolating}: for every $i$ and every $T\subset L\smin\partial L$ with $|T|\le i$,
the projection $A_i\to \prod_{x\in T} J^r_x(\RR^\ell,\RR^k)\cong P^T$ is onto.
\end{itemize}
\end{thm}

\medskip
As in section~\ref{sec-discrim} we define $B_i= A_i\smin \csigl$
and construct a ``resolution''
\begin{equation} \label{eqn-reso2} RB_i\to B_i~.
\end{equation}
In detail, $RB_i$ is the classifying space of
a topological poset whose elements are pairs $(f,T)$ where $f\in B_i$ and $T\subset L$ is a
bad event for $f$. The order relation is given by $(f,S)\le (g,T)$ iff $f=g$ and $S\subset T$.

\begin{lem} \label{lem-reso2} The space $RB_i$ is an ENR and the
resolution map $RB_i\to B_i$ is a proper homotopy equivalence.
\end{lem}

\proof It is shown in appendix~\ref{sec-badtop} that $RB_i$ is an ENR. Each fiber of $RB_i\to B_i$ is a
classifying space of a finite poset with maximal element, hence a contractible and compact simplicial complex.
This implies that $RB_i\to B_i$ is a \emph{cell-like}  map and therefore a proper homotopy equivalence.
\qed

\medskip
We filter $RB_i$ as follows. The underlying topological poset is filtered such that $(f,S)$ lives in the
$|S|$-th stage of the poset. This determines a filtration
\[ F_1 RB_i \subset F_2 RB_i \subset F_3 RB_i \subset \cdots \subset RB_i \]
such that all points in the interior of a nondegenerate simplex determined by a diagram
\[  (f,S_0) < (f,S_1) < \cdots < (f,S_{n-1}) < (f,S_n) \]
belong to $F_{|S_n|} RB_i$.

\begin{lem} \label{lem-resres2} Each $F_pRB_i$ is an ENR.
The resolution map $F_iRB_i \to B_i$ induces an isomorphism in locally finite homology
in dimensions $\ge d_i-\theta(i)$.
\end{lem}

\proof It is shown in in appendix~\ref{sec-badtop} that $F_pRB_i$ is an ENR. Fixing $i>0$,
let $C\subset B_i$ consist of all $f\in B_i$ which admit a bad event
of cardinality $\ge i$. Let $T=\{1,2,\dots,i\}$ and let
\[  Z\subset A_i\times\emb(T,L\smin\partial L) \]
consist of all $(f,e)$ such that $e(T)$ is a bad event for $f$.
Write $Z$ as a finite union of semi-algebraic subsets $Z_Q$ where $Q$ is a
sub-poset of the set of nonempty subsets of $T\cong e(T)$, and $Z_Q\subset Z$ consists of the
pairs $(f,e)$ such that the poset of bad sub-events of $e(T)$ is exactly $Q$.
Taking $j=\dim(BQ)$, we have that $Z_Q$ has codimension $\ge c(i,j)$ in
$A_i\times\emb(T,L\smin\partial L)$ by property (iv)
in theorem~\ref{thm-interpol2} (see~\ref{rem-codimestimate}), and
so has codimension $\ge c(i,j)-i\ell$ relative to $A_i$.
The set $Z_Q$ is a semialgebraic set
and it follows that its image in $A_i$ is also semialgebraic.
\newline
Now let $F_iRB_i|C$ be the portion of $F_iRB_i$ projecting to $C$.
Then $F_iRB_i|C$ is the image of a map
\[  \bigcup_Q~(Z_Q \times BQ) \lra F_iRB_i~. \]
The codimension of $Z_Q\times BQ$ relative to $A_i$ is $\ge c(i,j)-j-i\ell$ where $j=\dim(BQ)$.
Hence the codimension of $F_iRB_i|C$ relative to $A_i$ is $\ge \min_{j} (c(i,j)-i\ell-j)=\theta(i)$.
As the projection $F_iRB_i|C\to C$ is onto, it follows that the codimension of $C$ in $A_i$ is
$\ge \theta(i)$.
As in the proof of lemma~\ref{lem-resres}, it follows that
$H_*^\lf(F_iRB_i) \to  H_*^\lf(B_i)$
is an isomorphism when $*>d_i-\theta(i)$. \qed

\bigskip
There is a map
\begin{equation}\label{eqn-vassilspec2}
 F_p RB_i\smin F_{p -1}RB_i \lra \binom{L\smin\partial L}{p}
 \end{equation}
whose target is the space of unordered configurations of $p$ points in $L\smin\partial L$.

\begin{thm}\label{thm-loctriv2}
The map~(\ref{eqn-vassilspec2}) is a locally trivial projection, for $0<p\leq i$.
The fiber over a configuration $S$ is homeomorphic to the space of triples $(v,h,x)$ where
\begin{itemize}
\item[(a)] $v\in \RR^b$, with $b=d_i-p\dim(P)$;
\item[(b)] $h\in \prod_{s\in S} J^r_s(L,\RR^k)$ has $S$ as a bad event;
\item[(c)] $x$ belongs to the open cone on the classifying space of the poset of proper subsets of $S$
which are bad events for $h$.
\end{itemize}
\end{thm}

\proof
To show that it is locally trivial we factorize the map as follows:
\begin{equation} \label{eqn-vassilspecfac2}
F_p RB_i\smin F_{p-1}RB_i \lra E_p  \lra  \binom{L\smin\partial L}{p}.
\end{equation}
Here $E_p$ is the space of triples $(h,S,x)$ where $S$ is an unordered configuration of $p$ points
in $L\smin\partial L$,
\[ h\in \prod_{s\in S} J^r_s(L,\RR^k) \]
has $S$ as a bad event and $x$ is an element of the open cone mentioned in (c). The left-hand arrow associates
to $(f,S,x)$ as above the triple $(h,S,x)$, where $h=j_S^rf$. The right-hand arrow is forgetful, and it is
clearly a fiber bundle projection with fiber over $S$ homeomorphic to the space of pairs $(h,x)$ satisfying the conditions
(b) and (c). The left-hand arrow in~(\ref{eqn-vassilspecfac2})
is an affine space bundle with fibers $\RR^b$.
\qed

\begin{prop} \label{prop-speccompat2}
The spectral sequences in locally finite homology determined by the filtrations
\[  F_1 RB_i \subset F_2 RB_i \subset F_3 RB_i \subset \cdots \subset F_iRB_i \]
and
\[  F_1 RB_{i+1} \subset F_2 RB_{i+1} \subset F_3 RB_{i+1} \subset \cdots \subset F_iRB_{i+1} \]
are isomorphic up to a shift of $d_{i+1}-d_i$, equal to the dimension of $A_{i+1}/A_i$.
\end{prop}

\proof Analogous to the proof of proposition~\ref{prop-speccompat}. \qed

\begin{rem} \label{rem-belatedspec2}
{\rm There is a commutative ladder of homomorphisms
\begin{equation}\label{eqn-csiglladder2}
\xymatrix@R=12pt{
\vdots \ar[d]   &  \vdots  \\
\tilde H^*(A_{i+1}\cap\csigl)  \ar[d]   & H_{d_{i+1}-*-1}^\lf(F_{i+1}RB_{i+1}) \ar[l] \ar[u] \\
\tilde H^*(A_{i}\cap\csigl)  \ar[d]   & H_{d_{i}-*-1}^\lf(F_iRB_{i}) \ar[l] \ar[u] \\
\tilde H^*(A_{i-1}\cap\csigl)  \ar[d]   & H_{d_{i-1}-*-1}^\lf(F_{i-1}RB_{i-1}) \ar[l] \ar[u] \\
\vdots & \vdots \ar[u]
}
\end{equation}
The maps in the left-hand column are induced by the inclusions $A_i\to A_{i+1}$, and the maps
in the right-hand column are as in proposition~\ref{prop-speccompat2}. The horizontal map in row $i$ is the
composition
\[
\xymatrix{
H_{d_{i}-*-1}^\lf(F_iRB_{i}) \ar[r] & H_{d_{i}-*-1}^\lf(B_{i}) \ar[r]^-{\cong} & \tilde H^*(A_{i}\cap\csigl)
}
\]
By lemma~\ref{lem-resres2}, it is an isomorphism for $*+1<\theta(i)$. By the
dimension formula in theorem~\ref{thm-loctriv2}, the arrow
\[
H_{d_{i-1}-*-1}^\lf(F_{i-1}RB_{i-1}) \lra H_{d_{i}-*-1}^\lf(F_iRB_{i})
\]
is also an isomorphism when $*+1<\theta(i)$. Indeed, any deviation would be accounted for by
$F_iRB_i\smin F_{i-1}RB_i$, which by theorem~\ref{thm-loctriv2} has dimension
\[ \le i\ell+b+i\dim(P)-(\min_j c(i,j)-j)~=~i\ell+d_i-(\min_j c(i,j)-j)=d_i-\theta(i)~. \]
It follows that the arrow
$\tilde H^*(A_i\cap\csigl) \lra  \tilde H^*(A_{i-1}\cap\csigl)$
is an isomorphism for $*<\theta(i)-2$.
Finally, one concludes that the spectral sequence described in proposition~\ref{prop-speccompat2}
(stabilized with respect to $i$) \emph{converges} to the reduced cohomology of $\csigl$.
The dimension of $F_pRB_i\smin F_{p-1}RB_i$ is $\le d_i-\theta(p)$,
by theorem~\ref{thm-loctriv2} again. Therefore all $E^1_{p,q}$ vanish
where $p+q+d_i-1>d_i-\theta(p)$, which means
that $q\ge 2-p-\theta(p)$.
}
\end{rem}

\bigskip
We note that the analysis done of the spectral sequence converging to the reduced
cohomology of $\csigl(U)=C^\infty(U,\RR^k;\admi,\varphi)$
in section~\ref{sec-discman} generalizes painlessly to the more general setting of this
section where $\csigl(U)=C^\infty(U,\RR^k;\admib,\varphi)$.

\proof[Proof of~\ref{thm-vasscalcmultrel}]
It suffices to show that $\sS(U)\to T_p\sS(U)$ is approximately
$\theta(p)$-connected for every $U\in\sO\tm(L)$, where $\theta(p)=\min_j(c(p,j)-j-p\ell)$. We fix $p>0$ and some $i\gg p$. \newline
We ought to look at $\sS(U)$, the singular chain complex of $\csigl(U)$ as a contravariant
functor of $U$. Because of lemmas~\ref{lem-compactlim} and~\ref{lem-compactAdual} it is easier for us to
work for a little while with the reduced singular cochain complex of $\csigl(U)$ as a \emph{covariant}
functor of $U$. We have the following diagram of natural transformations between covariant functors
\[
\xymatrix@R=18pt{
{\tilde C^*(\csigl(U))} & &   {E_{i,p}(U)} \ar[d] \\
{\rule{0mm}{5mm}\hocolimsub{K} \tilde C^*(\csigl_K)} \ar[u]_-{\simeq} \ar[r]
& {\rule{0mm}{5mm}\hocolimsub{K}\tilde C^*(\csigl_K\cap A_i)}
\ar[r]^-{\simeq} & {E_{i,\infty}(U)}
}
\]
where $K$ runs through $\kappa(U)$.
We abbreviate this to
\[
\xymatrix{ Y & \ar[l]_-{\simeq} {Y^\sharp}  \ar[r] & {Y_i} \ar[r]^-{\simeq} & {E_{i,\infty}} & \ar[l] {E_{i,p}}
}
\]
suppressing the variable $U$. We apply the duality functor $D$ (see appendix~\ref{sec-polydual}) to obtain
\[
\xymatrix{ DY \ar[r]^-{\simeq} & {DY^\sharp} & \ar[l]{DY_i} & \ar[l]_-{\simeq}  {DE_{i,\infty}} \ar[r] &{DE_{i,p}}
}
\]
and form the Taylor approximations $T_p$~,
\[
\xymatrix@R=15pt{ DY \ar[r]^-{\simeq} \ar[d] & {DY^\sharp} \ar[d]
& \ar[l]{DY_i} \ar[d] & \ar[l]_-{\simeq}  {DE_{i,\infty}} \ar[r] \ar[d] &{DE_{i,p}} \ar[d]^-{\simeq} \\
T_pDY \ar[r]^-{\simeq} & {T_pDY^\sharp} & \ar[l]{T_pDY_i} &
\ar[l]_-{\simeq}  {T_pDE_{i,\infty}} \ar[r] &{T_pDE_{i,p}}~.
}
\]
In the top row, all maps are approximately $\theta(i)$-connected by remarks~\ref{rem-Kladder} and~\ref{rem-belatedspec},
except the one on the right which is approximately $\theta(p)$-connected by the same remarks.
In the bottom row, all maps except the one on the right are therefore approximately $(\theta(i)-\ell p-p)$-connected by lemma~\ref{lem-polydim}
in appendix~\ref{sec-polydual}. Noting that $T_pDE_{i,p}\simeq DE_{i,p}\simeq T_pDE_{i,i}$,
we can view the one on the right as being induced by the approximately $\theta(i)$-connected map
$DE_{i,\infty}\to DE_{i,i}$. It is therefore also approximately $(\theta(i)-\ell p-p)$-connected by lemma~\ref{lem-polydim}.
It follows that the left-hand
vertical arrow is approximately $\theta(p)$-connected. Since $DY$ is the reduced form of $\sS$, this completes the proof. \qed

\begin{rem} \label{rem-spec3} {\rm
As mentioned in the introduction, the fact that $\sS$ is analytic leads via general manifold calculus
principles to a spectral sequence converging to the
homology of $\sS(L)$. Of course we have already discussed that spectral sequence in~\ref{rem-belatedspec2},
but we have not given the standard manifold calculus description.
In order to use the description of the homogeneous layers
in~\cite{WeissEmb}, we replace $\sS$ by the spectrum-valued functor $F$ given by
\[  F(U) = \mathbf{H}\ZZ\,\wedge\!\csigl(U)_+ \]
for $U$ in $\sO\tm(L)$. For simplicity we also assume that $L$ is closed.
The $p$-th homogeneous layer of $F$ evaluated at $L$ is the spectrum of sections
with compact support of a fibered spectrum
\begin{equation} \label{eqn-Lambda} \mathbf{\Lambda}(p) \lra \binom{L}{p}~.
\end{equation}
The fiber of $\mathbf{\Lambda}(p)$ over a configuration $S\subset L$, where $|S|=p$, is
the total homotopy fiber of the $S$-cube
\[ T\mapsto F(N(T))\]
where $N(T)$ is a tubular neighborhood of $T$. In the category
of spectra, total homotopy fibers can be replaced by total homotopy cofibers up to a shift of $p$.
As $\csigl(N(\emptyset))$ is contractible and $\csigl(N(T))\simeq J^r_T(L,\RR^k;\admib)$ for nonempty $T\subset S$,
this means that the
fiber of $\mathbf{\Lambda}(p)$ over $S$ can be described as
\[ \Omega^p\mathbf{H}\ZZ~\wedge\!~\textup{cone}\big(\!\hocolimsub{\emptyset\ne T\subset S} J^r_T(L,\RR^k;\admib)\to \pt\big)~. \]
Here $J^r_T(L,\RR^k;\admib)\subset J^r_T(L,\RR^k)$ consists of the admissible elements; it is homeomorphic to
$\admi_T$. Finally using Poincar\'{e} duality,
the spectrum of sections of~(\ref{eqn-Lambda}) with compact support
can be identified with a twisted smash product with base equal to
that of~(\ref{eqn-Lambda}) and fibers obtained from those in~(\ref{eqn-Lambda}) by applying
$\Omega^W$, where $W$ is one of the tangent spaces of the configuration space. Its homotopy
group $\pi_{q-p}$ is therefore the homology group (twisted by the orientation character of the configuration space)
in degree $q-p+p+\ell p=q+\ell p\,\,$ of the projection map
\[      Y\lra \binom{L}{p} \]
whose fiber over a configuration $S$ is
\[  Y_S=\hocolimsub{\emptyset\ne T\subset S} J^r_T(L,\RR^k;\admib). \]  
(That homology group fits into a long exact sequence which also features the homology groups of $Y$ and the configuration space.)
This describes the term $E^1_{-p,q}$ of the spectral sequence. Clearly the emphasis here
is on admissible multijets, whereas the description in remark~\ref{rem-belatedspec2} has the
emphasis on bad multijets (bad events). 
}
\end{rem}

\begin{appendices}
\section{Interpolation and tameness}\label{sec-interpol}
In this Appendix we give proofs of theorem~\ref{thm-interpol} and
theorem~\ref{thm-interpol2}.

\medskip
Let $Q_n$ be the ring of polynomial maps $\RR^N\to \RR$ of degree $\le n$.

\begin{lem}\label{lem-idealtrans} For fixed $d>0$, there exists $n$ such that $Q_n$
is transverse to all codimension $d$ ideals of the space $C^\infty(\RR^N,\RR)$.
\end{lem}

\proof See for example Glaeser~\cite{Glaeser} or Vokrinek~\cite{Vokrinek}. The main point is that for every
codimension $d$ ideal $\mathfrak p$ there exist a finite set $S\subset \RR^N$ and a function
$v$ from $S$ to $\{1,\dots,d\}$ such that
\begin{equation}\label{eqn-idealdecomp}
 \prod_{x\in S} \mathfrak{m}_x^{v(x)} \subset \mathfrak p \subset \prod_{x\in S} \mathfrak{m}_x~.
\end{equation}
Here $\mathfrak{m}_x$ is the maximal ideal $\{f\in C^\infty(\RR^N,\RR)~|~f(x)=0\}$. \qed

\medskip
For a codimension $d$ ideal $\mathfrak p\subset  C^\infty(\RR^N,\RR)$, sandwiched as
in (\ref{eqn-idealdecomp}), the set $S$ will be called the \emph{spectrum} of $\mathfrak p$.

\begin{lem}\label{lem-idealopen} Let $n,\,m,\,d\in\NN$ be fixed. The space of $m$-dimensional affine subspaces of
$Q_n$ which are transverse to all codimension $d$ ideals of $C^\infty(\RR^N,\RR)$
with spectrum contained in $L$ is an open set in the space of $m$-dimensional
affine subspaces of $Q_n$.
\end{lem}

\proof By~\cite{Vokrinek}, the space of codimension $d$ ideals with spectrum contained in $L$ is a
compact Hausdorff space $X$, and there is a canonical vector bundle $E$ over $X$, with $d$-dimensional
fibers, so that the fiber over
$\mathfrak{p}\in X$ is equal to the quotient ring $C^\infty(\RR^N,\RR)/\mathfrak{p}$.
An $m$-dimensional affine subspace $A$ of $Q_n$ is transverse to all codimension $d$ ideals of
$C^\infty(\RR^N,\RR)$ with spectrum contained in $L$ if and only if the canonical vector bundle
map $X\times A\to E$ is fiberwise epimorphic. Hence the space of these affine subspaces is an open subset of
the space (product of $Q_n$ itself and a Grassmannian) of all $m$-dimensional affine
subspaces of $Q_n$ . \qed

\begin{lem}\label{lem-idealopen2} Let $n,\,m,\,d\in\NN$ be fixed. The space of $m$-dimensional affine subspaces of
$Q_n^k$ which are transverse to $\mathfrak p_1\times \mathfrak p_2\times\dots\times \mathfrak p_k$ for all
$k$-tuples $(\mathfrak p_1,\dots,\mathfrak p_k)$ of codimension $d$ ideals
of $C^\infty(\RR^N,\RR)$, all with spectrum contained in $L$, is an open set in the space of $m$-dimensional
affine subspaces of $Q_n^k$.
\end{lem}

\proof Same as for the lemma~\ref{lem-idealopen} above.
\qed

\bigskip
Define $g_L\co L\to \RR$ by $g_L=\sum_{i=1}^s f_i^r$~,
where $f_1,\dots,f_s$ are polynomial functions on $\RR^N$ which generate the ideal defining
$\partial L$ as a nonsingular algebraic subset of $\RR^N$.
As $r$ is even, we have $g_L^{-1}(0)=\partial L$. The $(r+1)$-th partial derivative of $g_L$
in the direction tangential to $L$ and normal to $\partial L$ is everywhere nonzero on $\partial L$.

\begin{lem} \label{lem-stilltransverse}
Let $i\in\NN$ and $d=i\dim(P)$. Suppose that $K\subset Q_n^k$ is a finite dimensional
linear subspace, transverse to $\mathfrak p_1\times\cdots \times\mathfrak p_k$ for all
tuples $(\mathfrak p_1,\dots,\mathfrak p_k)$ of codimension $\le d$ ideals
of $C^\infty(\RR^N,\RR)$, all with spectrum contained in $L$. Then $\varphi+g_L\cdot K \subset Sm$,
and for every subset $T$ of $L\smin\partial L$ with $|T|\le i$, the projection
\[  \varphi+g_L\cdot K \lra  \prod_{x\in T} J^r_x(L,\RR^k)~\cong~ P^T  \]
is onto.
\end{lem}

\proof For $T\subset L\smin\partial L$ with $|T|\le i$, let
$\mathfrak p_1=\mathfrak p_2=\cdots =\mathfrak p_k$ be
the ideal of $C^\infty(\RR^N,\RR)$ consisting of all functions $f$ whose
jet $j^r_xf$ vanishes for all $x\in T$. This has codimension $\le d$. By assumption,
$K$ is transverse to $\mathfrak p_1\times\cdots \times\mathfrak p_k$. It follows that
the projection
\[   K \lra  \prod_{x\in T} J^r_x(L,\RR^k)  \]
is onto. Since $g_L$ is nonzero at all points of $T$, we deduce that
\[  \varphi+g_L\cdot K \lra  \prod_{x\in T} J^r_x(L,\RR^k)  \]
is onto. \qed

\medskip
We will eventually construct the affine spaces $A_i$ of theorem~\ref{thm-interpol}
in the form $A_i=\varphi+g_L\cdot K_i$ for suitable $n\gg 0$ (depending on $i$) and finite dimensional
linear subspaces $K_i\subset Q_n^k$.
The three lemmas above show, broadly speaking, that the interpolation property~\ref{thm-interpol} (iv)
restricts our choice of $K_i$ to a \emph{non-empty open} collection of finite-dimensional linear spaces. In the following
we aim to show that the tameness property~\ref{thm-interpol} (iii) restricts our choice of $K_i$ to a \emph{dense} collection.

\medskip
Fix a linear space $K$ of polynomial maps $\RR^N\to \RR^k$, of finite dimension $m_0$.
Suppose that $K$ is \emph{tame} in the following sense: there exists
$\alpha\in\NN$ such that for every $f\in K$, the number of inadmissible
points of $\varphi+g_L\cdot f$ on $L$ is $<\alpha$. Fix such an $\alpha$.
(A large $\alpha$ is preferred; the reasons are given below.)
Next, choose $n$ sufficiently large so that the projection
\[ Q_n \to \prod_{x\in T} J^r_x(L,\RR^k) \]
is onto for every finite subset $T\subset L$ of cardinality at most $\alpha$,
and so that $K\subset Q_n^k$.
By lemma~\ref{lem-idealtrans}, such an $n$ exists.
We fix $m_1>m_0$ and turn our attention to the Grassmannian $Y$ of $m_1$-dimensional linear subspaces
of $Q_n^k$ containing $K$.
Let
\[  C_{K}\subset Y \]
be the semialgebraic subset consisting of all $K'$ which are \emph{wild} (as opposed to tame)
in the sense that there is $f\in K'$ such that $\varphi+g_L\cdot f$
has $\ge\alpha$ inadmissible points on $L$.

\begin{lem}\label{lem-famouscodim} If $2\alpha>m_1$, then
$C_{K}$ is nowhere dense in $Y$.
\end{lem}

\proof Let $EY$ be the tautological vector bundle on $Y$, with fiber $K'$ over the point $K'\in Y$.
Let $Z$ be the set of all $(K',f,S)$ where $K'\in Y$,
$f\in K'$ and $S\subset L$ has cardinality $\alpha$ and consists of inadmissible points for
$\varphi+g_L\cdot f$ (so that $S\cap \partial L=\emptyset$).
Then
\[  Z~\subset~EY \times L^\alpha \]
as a semialgebraic subset.
By the interpolation property which we are assuming for
$Q_n$~, and by our assumption on $\admi$, the codimension of this semialgebraic subset is at least
$\alpha\cdot(\ell+2)$ (see remark~\ref{rem-codimestimate}).
The image of $Z$ under the projection to $Y$ therefore still has codimension at least
$2\alpha-m_1$ everywhere. That image is precisely $C_{K}$.
Therefore the codimension of $C_{K}$ in $Y$ is at least $1$ if
$2\alpha>m_1$.
\qed

\bigskip
Now we construct the finite dimensional linear spaces $K_i$~, and the affine spaces
$A_i=\varphi+g_L\cdot K_i$ of theorem~\ref{thm-interpol}. Suppose per induction that
$K_i$ has already been constructed for a specific $i\ge 0$, has dimension $m_0$
and consists of polynomial maps from $\RR^N$ to $\RR^k$. We also assume that
the following condition $(\lambda_i)$ is satisfied: for
every $k$-tuple $f=(f_1,\dots,f_k)$ of monomials of degree $\le i$ in $N$ variables, there exists a map in $K_i$
whose distance from $f$ in the $C^i$ topology (on the space of smooth maps from $L$ to $\RR^k$)
is less than $2^{-i}$. 
\newline
Choose $s\in \NN$ large enough so that $K_i$ is properly contained in $Q_s^k$ and
$Q_s$ is transverse to all codimension $\le d$ ideals of the
space $C^\infty(\RR^N,\RR)$, where
\[ d=(i+1)\dim(P)~. \]
Let $m_1=\dim(Q_s^k)$.
Choose $\alpha>m_1/2$ and choose $n$ sufficiently large so that the projection
\[ Q_n\to \prod_{x\in T} J^r_x(L,\RR^k) \]
is onto for every finite subset $T\subset L$ of cardinality at most $\alpha$. Also, arrange $n>s$
and $s>i$.
By lemmas~\ref{lem-idealopen} and~\ref{lem-famouscodim}, there exists a linear subspace $K_{i+1}$
of $Q_n^k$ of dimension $m_1$ and as close to $Q_s^k$ as we might wish, with the following properties:
\begin{itemize}
\item[--] $K_i\subset K_{i+1}$;
\item[--] $K_{i+1}$ is transverse  to $\mathfrak p_1\times\dots\times \mathfrak p_k$ for all
tuples $(\mathfrak p_1,\dots,\mathfrak p_k)$ of
codimension $\le d$ ideals of $C^\infty(\RR^N,\RR)$, all with spectrum contained in $L$,
where
\[ d=(i+1)\dim(P)~; \]
\item[--] any set of inadmissible points, alias bad event, (in $L$) for any $f\in A_{i+1}:= \varphi+g_L\cdot K_{i+1}$ has cardinality $<\alpha$.
\end{itemize}
Then the tameness condition for
$A_{i+1}=\varphi+g_L\cdot K_{i+1}$ is satisfied. The interpolation condition
is also satisfied by lemma~\ref{lem-stilltransverse}. Since $K_{i+1}$ can be as close
as we wish to $Q_s$~, and $s>i$, we can choose it so that condition $(\lambda_{i+1})$
is satisfied. This completes the induction.
Because condition $(\lambda_i)$
is satisfied by $K_i$~, for all $i$, the union $\bigcup_i A_i$ is dense in $Sm$.

\begin{rem} \label{rem-codimestimate} {\rm
Suppose that $A_i$ is a finite dimensional affine subspace of $\Sm$ which has the interpolation property (iv) of theorem~\ref{thm-interpol}.
Let $T=\{1,2,\dots,i\}$ and let
\[  Z\subset A_i\times\emb(T,L\smin\partial L) \]
consist of all $(f,e)$ such that $e(T)$ is a bad event for $f$. Suppose we fix a point $(f,e)\in Z$ and
choose a smooth trivialization of the tangent bundle $TL$ near $e(T)$, not necessarily algebraic.
We can then use this to trivialize the jet bundle near $e(T)$ and construct a germ of maps
\[ A_i\times\emb(T,L\smin \partial L) \lra \prod_{x\in T}J_x^r(\RR^\ell,\RR^k)\cong P^T\]
(near $(f,e)$) by mapping an element $(g,b)$ to the multijet of $g$ at $b(T)$.
Then the germ of $Z$ near $(f,e)$ is the preimage of
$(P\smin \admi)^T$ under this map germ. On the other hand, by the interpolation property, the map germ is a submersion.
Therefore the codimension of the germ  of $Z$ near $(f,e)$, as a germ of manifold stratified subsets
of $A_i\times\emb(T,L\smin\partial L)$,
 is $\ge$ the codimension of $(P\smin \admi)^T$ in $P^T$, which is $\ge i(\ell+2)$.\newline
Similarly, in the situation of lemma~\ref{lem-famouscodim}, we have
\[  Z~\subset~EY \times L^\alpha~. \]
Given $(K',f,S)\in Z$, choose a trivialization of $TL$ near $f(S)$ and, as above, use this to construct a submersion germ
\[ (EY\times L^\alpha,(K',f,S)) \lra \prod_{x\in S}J_x^r(\RR^\ell,\RR^k)\cong P^S~\]
such that the germ of $Z$ near $(K',f,S)$ is the preimage of $(P\smin \admi)^S$ under this map germ.
Therefore the codimension of the germ of $Z$ near $(K',f,S)$ is $\ge \alpha (\ell+2)$.
}
\end{rem}

\medskip
Now we prove theorem~\ref{thm-interpol2}.
We will eventually construct the affine spaces $A_i$ of theorem~\ref{thm-interpol}
in the form $A_i=\varphi+g_L\cdot K_i$ for suitable $n\gg 0$ (depending on $i$) and finite dimensional
linear subspaces $K_i\subset Q_n^k$.
The three lemmas above show, broadly speaking, that the interpolation property~\ref{thm-interpol} (iv)
restricts our choice of $K_i$ to a \emph{non-empty open} collection of finite-dimensional linear spaces. In the following
we aim to show that the tameness property~\ref{thm-interpol} (iii) restricts our choice of $K_i$ to a \emph{dense} collection.

\medskip
Fix a linear space $K$ of polynomial maps $\RR^N\to \RR^k$, of finite dimension $m_0$.
Suppose that $K$ is \emph{tame} in the following sense: there exists
$\alpha\in\NN$ such that for every $f\in K$, the cardinality of any bad event for $\varphi+g_L\cdot f$ on $L$ is $<\alpha$.
Fix such an $\alpha$. Next, choose $n$ sufficiently large so that the projection
\[ Q_n \to \prod_{x\in T} J^r_x(L,\RR^k) \]
is onto for every finite subset $T\subset L$ of cardinality at most $\alpha+\zeta-1$ where $\zeta$ is the upper bound of definition~\ref{defn-codim},
and so that $K\subset Q_n^k$.
By lemma~\ref{lem-idealtrans}, such an $n$ exists.
We fix $m_1>m_0$ and turn our attention to the Grassmannian $Y$ of $m_1$-dimensional linear subspaces
of $Q_n^k$ containing $K$.
Let
\[  C_{K}\subset Y \]
be the semialgebraic subset consisting of all $K'$ which are \emph{wild} in the sense that there is $f\in K'$ such that $\varphi+g_L\cdot f$
has a bad event of cardinality $\ge\alpha$ on $L$.

\begin{lem}\label{lem-famouscodim2} If $2\alpha>m_1$, then
$C_{K}$ is nowhere dense in $Y$.
\end{lem}

\proof Let $EY$ be the tautological vector bundle on $Y$, with fiber $K'$ over the point $K'\in Y$.
Let $Z$ be the set of all $(K',f,S)$ where $K'\in Y$ and
$f\in K'$ and $S\subset L$ a bad event for   
$\varphi+g_L\cdot f$ whose cardinality is between $\alpha$ and $\alpha+\zeta-1$.
Then
\[  Z~\subset~EY \times L^\alpha \]
as a semialgebraic subset. From here onwards we proceed as in the proof of lemma~\ref{lem-famouscodim}.
\qed

\bigskip
Now we construct the finite dimensional linear spaces $K_i$~, and the affine spaces
$A_i=\varphi+g_L\cdot K_i$ of theorem~\ref{thm-interpol2}. Suppose per induction that
$K_i$ has already been constructed for a specific $i\ge 0$, has dimension $m_0$
and consists of polynomial maps from $\RR^N$ to $\RR^k$. We also assume that
the following condition $(\lambda_i)$ is satisfied: for
every $k$-tuple $f=(f_1,\dots,f_k)$ of monomials of degree $\le i$ in $N$ variables, there exists a map in $K_i$
whose distance from $f$ in the $C^i$ topology (on the space of smooth maps from $L$ to $\RR^k$)
is less than $2^{-i}$. 
\newline
Choose $s\in \NN$ large enough so that $K_i$ is properly contained in $Q_s^k$ and
$Q_s$ is transverse to all codimension $\le d$ ideals of the
space $C^\infty(\RR^N,\RR)$, where
\[ d=(i+1)\dim(P)~. \]
Let $m_1=\dim(Q_s^k)$.
Choose $\alpha>m_1/2$ and choose $n$ sufficiently large so that the projection
\[ Q_n\to \prod_{x\in T} J^r_x(L,\RR^k) \]
is onto for every finite subset $T\subset L$ of cardinality at most $\alpha+\zeta-1$. Also, arrange $n>s$
and $s>i$.
By lemmas~\ref{lem-idealopen} and~\ref{lem-famouscodim2}, there exists a linear subspace $K_{i+1}$
of $Q_n^k$ of dimension $m_1$ and as close to $Q_s^k$ as we might wish, with the following properties:
\begin{itemize}
\item[--] $K_i\subset K_{i+1}$;
\item[--] $K_{i+1}$ is transverse  to $\mathfrak p_1\times\dots\times \mathfrak p_k$ for all
tuples $(\mathfrak p_1,\dots,\mathfrak p_k)$ of
codimension $\le d$ ideals of $C^\infty(\RR^N,\RR)$, all with spectrum contained in $L$,
where
\[ d=(i+1)\dim(P)~; \]
\item[--] any bad event (in $L$) for any $f\in A_{i+1}:= \varphi+g_L\cdot K_{i+1}$ has cardinality $<\alpha$.
\end{itemize}
Then the tameness condition for
$A_{i+1}=\varphi+g_L\cdot K_{i+1}$ is satisfied. The interpolation condition
is also satisfied by lemma~\ref{lem-stilltransverse}. Since $K_{i+1}$ can be as close
as we wish to $Q_s$~, and $s>i$, we can choose it so that condition $(\lambda_{i+1})$
is satisfied. This completes the induction.
Because condition $(\lambda_i)$
is satisfied by $K_i$~, for all $i$, the union $\bigcup_i A_i$ is dense in $Sm$.

\section{Topology of badness} \label{sec-badtop}
Fix a positive integer $\alpha$ and let $X$ be the space of nonempty subsets of $\RR^N$ of cardinality at most $\alpha$,
equipped with the Hausdorff metric. Let $Y$ be the space of all maps from
$\{1,2,3,\dots,\alpha\}$ to $\RR^N$, topologized as the $\alpha$-fold power of $\RR^N$. There is a surjective map
$q\co Y\lra X$ taking $f\in Y$ to its image $\im(f)$.

\begin{lem} \label{lem-Haus} The map $q$ is an identification map.
\end{lem}

\proof For fixed $f\in Y$ and
$\varepsilon>0$, the set of $g\in Y$ whose image $\im(g)$ has Hausdorff
distance $<\varepsilon$ from $\im(f)$
is clearly open in $Y$. This means that $q$ is continuous.  \newline
Next, let $S\in X$ be fixed and let $U\subset X$ be such that
$S\in U$. Suppose that $U$ is not a neighborhood of $S$. We need to show that $q^{-1}(U)$ is not a
neighborhood of $q^{-1}(S)$. Since $U$ is not a neighborhood of $S$, there exists a
sequence $(S_i)_{i\in\NN}$ in $X$ converging
to $S$ such that $S_i\notin U$ for all $i$. Choose $\varepsilon>0$ such that the $\varepsilon$-neighborhoods
of the $s\in S$ are pairwise disjoint in $\RR^N$. Each $S_i$ determines a function $m_i\co S\to \NN$,
where $m_i(s)$ is the number of elements of $S_i$ which have distance $<\varepsilon$
from $s\in S$. As $\sum_{t\in S} m_i(t)\le \alpha$, there are only finitely many possibilities for the $m_i$. Hence the
sequence $(S_i)_{i\in\NN}$ has a subsequence $(S_{i_r})_{r\in\NN}$ such that $m_{i_r}$ is the same for all $r$.
We rename: $T_r=S_{i_r}$. Now it is easy to lift the sequence $(T_r)_{r\in\NN}$ to a convergent sequence $(f_r)_{r\in\NN}$
in $Y$, with a limit which we call $f$. Then $f\in q^{-1}(S)$ but
$f_r\notin q^{-1}(U)$, so that $q^{-1}(U)$ is not a neighborhood of $q^{-1}(S)$. \qed

\bigskip
The equivalence relation on $Y$ which has $f\sim g$ if and only if $\im(f)=\im(g)$ is fairly
easy to understand. In particular the following is easy.

\begin{lem} \label{lem-triang} There is a triangulation of $Y$ which agrees with the linear structure {\rm(}every
simplex is the convex hull in $Y$ of its vertex set{\rm)} and which descends to a
triangulation of the quotient $Y/\sim$~.
\end{lem}

\medskip
Next let $A$ be a finite dimensional affine subspace of $\Sm$ (as in section~\ref{sec-discrim}
or section~\ref{sec-multdiscrim}) which consists of polynomial maps.
Assume that $A$ is tame in the sense that the cardinality of any bad event of any $f\in A$ is
not greater than $\alpha$. Let $B\subset A$ be the subset of all $f\in A$ which have a bad event
(see section~\ref{sec-discrim} or definition~\ref{defn-mapbadness}).
Let $\Lambda$ be the poset of all pairs $(f,S)$ such that $f\in B$ and $S\subset L\smin \partial L$ is a bad
event for $f$. The order relation has $(f,S)\le (g,T)$ if and only if $f=g$ and $S\subset T$. We can write
\[  \Lambda \subset B\times X~. \]

\begin{cor} $\Lambda$ is an ENR and the projection $\Lambda \to B$ is proper.
\end{cor}

\proof By \cite[Ch.~VI, thm.~1.2]{Hu}, to show that $\Lambda$ is an ENR and closed in $B\times X$,
it is enough to show that $\Lambda \cap (B\times Q)$ is a closed ENR in $B\times Q$
for every simplex $Q\subset X$. By lemma~\ref{lem-triang} we can lift $Q$ to a simplex
$\bar Q$ in the vector space $Y$. It remains to show that the preimage of $\Lambda\subset B\times X$ under
the composite map
\[
\xymatrix{  B\times \bar Q \ar[r]^{\text{incl.}} & B\times Y \ar[r]^{\id\times q} & B\times X
}
\]
is an ENR. That preimage is clearly a closed semialgebraic subset of $A\times Y$. \newline
For the properness statement we introduce $Z$~, the space of subsets of $L$ of cardinality
$\le \alpha$. This is a quotient of $L^\alpha$. Hence it is compact.
The inclusion of $RB$ in $B\times X$ factors as
\[  RB \lra B\times Z \lra B\times X~. \]
By the above, $RB$ is closed in $B\times Z$ since it is closed in $B\times X$.
The projection $B\times Z\to B$ is proper since $Z$ is compact. Therefore
its restriction to $RB$ is also proper. \qed

\begin{cor} For $s\ge 0$, let $N_s\Lambda$ be the set of all order-preserving maps from $\{0,1,\dots,s\}$ to $\Lambda$.
Then $N_s\Lambda$ is an ENR and the projection $N_s\Lambda\to B$ is proper.
\end{cor}

\proof This is similar to the previous proof. We do the case $s=1$. Let $Q_0$ and $Q_1$ be simplices in $X$
with lifts $\bar Q_0$ and $\bar Q_1$ to $Y$. We need to show that the preimage in
\[ (B\times \bar Q_0)\times (B\times \bar Q_1) \]
of the order relation in $\Lambda\times \Lambda\subset (B\times X)\times(B\times X)$
is a closed sub-ENR. That preimage is
a closed semialgebraic subset of $(A\times Y)\times(A\times Y)$. Hence it is an ENR. We conclude as
before that $N_1\Lambda$ is an ENR, closed in $\Lambda\times\Lambda$. The inclusion
of $N_1\Lambda$ in $\Lambda\times\Lambda$ factors as
\[  N_1\Lambda \lra B\times Z \times Z \lra \Lambda\times \Lambda~, \]
with $Z$ as in the previous proof.
Hence $N_1\Lambda$ is closed in $B\times Z\times Z$. Since the projection from
$B\times Z\times Z$ to $B$ is proper, its restriction to $N_1\Lambda$ is also proper.

\qed

\bigskip
Finally we consider $\Lambda$ as a topological poset, with a classifying space which we call
$RB$ to keep the analogy with section~\ref{sec-discrim}. It is the geometric realization
of a simplicial space $s\mapsto N_s\Lambda$. Each $N_s\Lambda$ is an ENR and
all simplicial operators $N_s\Lambda\to N_t\Lambda$ are proper.
Because the length of any chain of non-identity morphisms in $\Lambda$ is bounded by $\alpha$,
every element of $N_s\Lambda$ for $s>\alpha$ is in the image of some degeneracy operator.
These facts together with \cite[Ch. VI, thm~1.2]{Hu} imply the following.

\begin{prop} The classifying space $RB$ of $\Lambda$ is an ENR and the projection
$RB\to B$ is proper.
\end{prop}

\section{Excision in locally finite homology} \label{sec-exci}
Let $A$ be an ENR, and $B$ a closed subset of $A$ which is also an ENR.

\begin{prop} \label{prop-lfex} There is an isomorphism $H_*^\lf(A,B) \cong H_*^\lf(A\smin B)$.
\end{prop}

\proof This is well known. It can be proved using (unreduced) Steenrod homology, a homology theory $H_*^\ste$ for
pairs of compact metric spaces \cite{Ste1,Ste2,Milnor61,Ferry}. Milnor showed \cite{Milnor61} that $H_*^\ste$ satisfies,
in addition to the seven Eilenberg-Steenrod axioms for a homology theory, two further axioms numbered eight and nine.
Axiom eight is a wedge axiom
for a type of infinite wedge, while axiom nine states that for a pair $(X,Y)$ of compact metric spaces,
the identification map $(X,Y)\to (X/Y,\pt)$ induces an isomorphism in Steenrod homology. These axioms
imply that for a pair of finite dimensional locally finite CW-spaces $(A,B)$ with one-point
compactifications $A^\omega$ and $B^\omega$,
the Steenrod homology of the pair $(A^\omega,B^\omega)$ is isomorphic to the locally finite homology of
the pair $(A,B)$. Indeed the spectral sequence in Steenrod homology determined by the skeleton filtration
of $A$ relative to $B$ collapses (by the axioms) at the $E^2$ page, which gives the isomorphism. It follows
that $H_*^\ste(A^\omega,B^\omega)$ is also
isomorphic to the locally finite homology of $(A,B)$ when $(A,B)$ is a pair of ENRs with $B$ closed in $A$.
But by axiom number nine, $H_*^\ste(A^\omega,B^\omega)$ is also isomorphic to $H_*^\ste(A/B,\pt)$ which
in turn is isomorphic to the locally finite homology of $A\smin B$. \qed

\begin{cor} \label{cor-Adual} Let $B$ be a closed ENR in $\RR^d$. Then
\[  \tilde H^*(\RR^d\smin B)\cong H_{d-*-1}^\lf(B) \]
where $\tilde H^*$ is reduced cohomology.
\end{cor}

\proof The long exact sequence in locally finite homology of the pair $(\RR^d,B)$ shows
\[ H_{d-s-1}^\lf(B)\cong H_{d-s}^\lf(\RR^d,B) \]
for $s>0$, because $H_*^\lf(\RR^d)$ is isomorphic
to $H^{d-*}(\RR^d)$ by Poincar\'{e} duality for noncompact manifolds. By proposition~\ref{prop-lfex}
and another application of Poincar\'{e} duality, we have
\[ H_{d-s}^\lf(\RR^d,B)\cong H_{d-s}^\lf(\RR^d\smin B)\cong H^s(\RR^d\smin B)~. \]
This isomorphism extends to a map of long exact sequences
\[
\xymatrix@R=11pt{  & \\
H_{d-s-1}^\lf(B) \ar[r] \ar@{..>}[u] &  H^{s+1}(\RR^d,\RR^d\smin B) \ar@{..>}[u] \\
H_{d-s}^\lf(\RR^d,B) \ar[r] \ar[u] &  H^s(\RR^d\smin B) \ar[u] \\
H_{d-s}^\lf(\RR^d) \ar[r] \ar[u] & H^s(\RR^d) \ar[u] \\
H_{d-s}^\lf(B) \ar[r] \ar[u] & H^s(\RR^d,\RR^d\smin B) \ar[u] \\
 \ar@{..>}[u]  &  \ar@{..>}[u]
}
\]
which, by the five lemma, is also an isomorphism. In particular
\[\qquad\qquad\qquad
H_{d-s-1}^\lf(B)\cong H^{s+1}(\RR^d,\RR^d\smin B)\cong H^s(\RR^d\smin B)~.\qquad\qquad\qquad \qed \]

\begin{expl} {\rm Something like the ENR condition in corollary~\ref{cor-Adual} is necessary. The set
$B=\{n^{-1}\in \RR~|~n=1,2,3,\dots\}\cup\{0\}$
is closed in $\RR$, but it is not an ENR. It is compact, and so
\[ H^\lf_0(B)=H_0(B)\cong \bigoplus_{s\in B} \ZZ~. \]
By contrast, $H^0(\RR\smin B)\cong \prod_{s\in B} \ZZ$. The two groups are not isomorphic since one is uncountable and
the other is countable.
}
\end{expl}

\begin{rem} \label{rem-subd} {\rm
There are situations where we need to promote proposition~\ref{prop-lfex} to a statement at the chain level.
Therefore we shall briefly describe a chain map of locally finite singular chain complexes
\begin{equation} \label{eqn-lcsubdivide}
 C_*^\lf(A,B) \lra C_*^\lf(A\smin B),
\end{equation}
defined when $(A,B)$ is a pair of locally compact spaces, with $B$ closed in $A$. By
$C_*^\lf(A,B)$ we mean the quotient $C_*^\lf(A)/C_*^\lf(B)$, and we note that the
inclusion $C_*^\lf(B)\to C_*^\lf(A)$ is a cofibration, i.e., split injective in each degree.
If $A$ and $B$ are ENRs, then~(\ref{eqn-lcsubdivide}) induces an isomorphism in homology, the isomorphism of
proposition~\ref{prop-lfex}. \newline
The idea is simple.
Given any singular simplex $f\co \Delta^q\to A$, we have the closed set $f^{-1}(B)\subset \Delta^q$. For $p\ge 0$ let
$S_p$ be the set of $q$-simplices in the $p$-th barycentric subdivision of $\Delta^q$ which have empty intersection
with $f^{-1}(B)$ and which are \emph{not} contained in a simplex of the $(p-1)$-th subdivision which also has
empty intersection with $f^{-1}(B)$. Then the formal sum
\[  \sum_{p\ge 0\,}\sum_{\sigma\in S_p} \sigma \]
is a locally finite $q$-chain in $A\smin B$. By mapping the $q$-simplex $f$ to that locally finite $q$-chain,
we obtain a chain map of the form~(\ref{eqn-lcsubdivide}). It has the following naturality property: Given
locally compact $A$ and closed subsets $B_0$ and $B_1$ with $B_0\subset B_1$, there is a commutative diagram
\[
\xymatrix{
C_*^\lf(A,B_0) \ar[r] \ar[d] & C_*^\lf(A\smin B_0) \ar[d] \\
C_*^\lf(A,B_1) \ar[r] & C_*^\lf(A\smin B_1)
}
\]
where the horizontal arrows are of the form~(\ref{eqn-lcsubdivide}) and the right-hand vertical arrow
is also of that type (note that $A\smin B_1=(A\smin B_0)\smin(B_1\smin B_0)$). \newline
In the case where $A=\RR^d$, the chain map~(\ref{eqn-lcsubdivide}) fits (by naturality) into a commutative
diagram of chain maps
\[
\xymatrix{
C_*^\lf(\RR^d) \ar[r]^-{=}  \ar[d] & C_*^\lf(\RR^d) \ar[d] & \ar[l]_-{\simeq} C^{d-*}(\RR^d) \ar[d]  \\
C_*^\lf(\RR^d,B) \ar[r]^-{\simeq} & C_*^\lf(\RR^d\smin B) & \ar[l]_-{\simeq} C^{d-*}(\RR^d\smin B)
}
\]
Passing to the mapping cones of the vertical arrows, we obtain a chain of natural homotopy equivalences
\[  C_*^\lf(B) \lra \cdots \longleftarrow \tilde C^{d-*-1}(\RR^d\smin B) \]
where the tilde indicates a reduced cochain complex. (The reduced singular cochain complex of a space
$X$ is the mapping cone of the chain map $C^*(\pt)\to C^*(X)$ induced by the unique map
$X\to\pt$.)
}
\end{rem}

We state one more result about locally finite homology which does not depend on ENR assumptions.

\begin{prop} \label{prop-flasque} Let $Y$ be any locally compact metrizable space such that
the one-point compactification is also metrizable. Then $H_*^\lf(Y\times(0,1]\,)$ is zero.
\end{prop}

\proof Let $X$ be the one-point compactification of $Y\times[0,1]$ and let $X'$ be the
one-point compactification of $Y\times 0$. By axiom $8^\lf$ in \cite{Milnor61} for the Steenrod homology groups,
\[  H_q^\ste(X,X') \cong H_q^\lf(X\smin X')\cong H_q^\lf(Y\times(0,1]\,)~. \]
Therefore it is enough to show that $H_q^\ste(X,X')=0$. But this follows from the homotopy invariance property
of the Steenrod homology groups (which is part of the seven Eilenberg-Steenrod axioms). \qed

\section{Boundary conditions} \label{sec-boundary}
Let $\varphi$ be an admissible smooth map from $L$ to $\RR^k$. Here \emph{admissible} means
$\admi_\bullet$-admissible as in section~\ref{sec-multdiscrim}; we recall that this generalizes $\admi$-admissible
as in section~\ref{sec-discrim}. In our definition of the space $C^\infty(U,\RR^k;\admi_\bullet,\varphi)$
there is a boundary condition (a): a smooth map $f\co L\to \RR^k$ which is everywhere admissible belongs to that
space if its $r$-jet at every $x\in \partial L$ agrees with the $r$-jet of $\varphi$ at $x$. This deviates
slightly from the original definition \cite{Vassiliev1,Vassiliev2}. Vassiliev has the condition (b) that
the full Taylor series of $f$ at every $x\in \partial L$ agree with the full Taylor series of $\varphi$ at $x$
(in local coordinates at $x$). Condition (b) has some technical advantages, as we shall see
in the proof of proposition~\ref{prop-bdryindep} below. Our weaker condition (a) also has some technical advantages.
It is much easier to find polynomial maps which satisfy it, in the case where $L\subset \RR^k$ is a semi-algebraic
subset (details as in section~\ref{sec-discrim}).

\begin{lem} \label{lem-bdrystrgwk} The two types of boundary conditions, (a) and (b), lead to
weakly homotopy equivalent versions of $C^\infty(U,\RR^k;\admi_\bullet,\varphi)$.
\end{lem}

\proof Denote the two versions by $X_a$ and $X_b$~, respectively. The inclusion $X_b\to X_a$ is claimed
to be a weak homotopy equivalence. To show this, let $K\subset X_a$ be any compact subset (for example,
the image of a continuous map from a compact CW-space to $X_a$). Choose a compact collar $Q$
on $\partial L$ and identify if with $\partial L\times [0,1]$, relative to $\partial L\cong L\times 0$.
Choose a smooth function $\psi\co [0,1]\to [0,1]$ so that $\psi(t)=0$ for $t$ close to $0$ and
$\psi(t)=1$ for $t$ close to $1$. Define $\Psi\co L\to [0,1]$ so that $\Psi(x,t)=\psi(t)$
for $(x,t)$ in $\partial L\times[0,1]\cong Q\subset L$ and
$\Psi(y)=1$ for all $y\in L\smin Q$.
We ask whether, for $s\in [0,1]$ and $f\in K$, the map
\[  f_s:=(1-s)f+s(\Psi\cdot f+(1-\Psi)\varphi) \]
from $L$ to $\RR^k$ is everywhere admissible. If so, then $f\mapsto (f_s)_{s\in[0,1]}$ defines a homotopy from the inclusion
$K\to X_a$ to a map which lands in the subspace $X_b$. Moreover, if $K$ is already contained in $X_b$,
then the homotopy would not take it out of $X_b$. In short, if the formula for $f_s$ is good, then the
proof is complete. \newline
Of course there is no guarantee that the formula is good as it stands. But we can easily improve on it.
Let $Q^\varepsilon=\partial L\times[0,\varepsilon]\subset \partial L\times[0,1]\cong Q$ and define $\Psi^\varepsilon\co L\to [0,1]$
by $\Psi^\varepsilon(x,t)=\Psi(x,t/\varepsilon)$ for $(x,t)\in Q^\varepsilon$ and $\Psi^\varepsilon(y)=1$
for all $y\in L\smin Q^\varepsilon$. A calculation which we leave to the reader shows that
the above formula for $f_s$ turns into a good formula if we replace $\Psi$ by $\Psi^\varepsilon$ for sufficiently small
$\varepsilon$. This is based on the existence of a constant $c>0$ such that for $f\in K$ and $p$ with
$0\le p\le r$, all $p$-th partial derivatives
of $f-\varphi$ at a point $(x,t)\in \partial L\times[0,1]\cong Q\subset L$ are bounded
in size by $ct^{r+1-p}$. \qed

\begin{rem} {\rm Lemma~\ref{lem-bdrystrgwk} justifies a statement made in section~\ref{sec-discrim}, namely,
there is no loss of generality in assuming that $r$ is even. We recall that $\admi_T\subset P^T$ for a finite
nonempty set $T$, where $P$ is the vector space of polynomial maps of degree $\le r$ from $\RR^\ell$ to $\RR^k$.
Let $\mathfrak Y^T$
be the preimage of $\admi^T$ under the projection $P_1^T\to P^T$, where $P_1$ is the vector space of
polynomial maps of degree $\le r+1$ from $\RR^\ell$ to $\RR^k$. There is an inclusion
\[  C^\infty(U,\RR^k;\mathfrak Y_\bullet,\varphi) \lra C^\infty(U,\RR^k;\admi_\bullet,\varphi)~. \]
We are using type (a) boundary conditions in the sense of lemma~\ref{lem-bdrystrgwk}.
If $L$ has empty boundary, this inclusion is an identity. If $\partial L$ is nonempty, it would still be an identity
if type (b) boundary conditions were used; therefore it is a weak equivalence with type (a) boundary conditions.
}
\end{rem}

\bigskip
Next, let $\varphi_0$ and $\varphi_1$ be two admissible smooth maps from $L$ to $\RR^k$. Here \emph{admissible} means
$\admi_\bullet$-admissible as in section~\ref{sec-multdiscrim}; we recall that this generalizes $\admi$-admissible
as in section~\ref{sec-discrim}.

\begin{prop} \label{prop-bdryindep} The contravariant functors $F_0$ and $F_1$ on $\sO\tm(L)$ defined by
\[  F_0(U)=C^\infty(U,\RR^k;\admi_\bullet,\varphi_0)~, \qquad
F_1(U)= C^\infty(U,\RR^k;\admi_\bullet,\varphi_1) \]
are related by a chain of natural weak homotopy equivalences.
\end{prop}

The proof is an elaboration of an argument given in \cite{Vassiliev2}.
As a preparation, we choose a compact collar $Q$ on $\partial L$, so that $Q\subset L$ and
$\partial L\subset Q\cong \partial L\times[0,2]$. Let $\sO\tm(L;Q)\subset \sO\tm(L)$ be the full sub-poset
consisting of the $U\in \sO\tm(L)$ which contain the collar $Q$.

\begin{lem} \label{lem-bdryindep}
There is a natural weak homotopy equivalence from $F_1|\sO\tm(L;Q)$ to $F_0|\sO\tm(L;Q)$.
\end{lem}

\proof We use Vassiliev's boundary condition (b) throughout, as in lemma~\ref{lem-bdrystrgwk}.
The collar $Q$ can be identified with $\partial L\times[0,2]$, relative to
$\partial L\cong \partial L\times 0$. Choose $\varepsilon>0$ and a smooth embedding $\psi\co [0,2]\to[0,2]$ which extends
\begin{itemize}
\item the map $t\mapsto t+1$ on $[0,2\varepsilon]$;
\item the map $t\mapsto t+\frac{1}{2}$ on $[1,1+2\varepsilon]$;
\item the identity on $[2-\varepsilon,2]$.
\end{itemize}
Choose a smooth admissible map $\varphi_2$ from $\partial L\times[0,1]$ to
$\RR^k$ such that $\varphi_2(x,t)=\varphi_0(x,t)$ for $t\in[0,\varepsilon]$ and
such that the map from $\partial L\times[0,2]$ to $\RR^k$ defined by $(x,t)\mapsto \varphi_2(x,t)$
for $t\in[0,1]$ and $(x,t)\mapsto\varphi_1(x,t-1)$ for $t\in[1,2]$ is smooth. \newline
Given $U$ in $\sO\tm(L)$ containing $Q$,
and $f\in F_1(U)$, define $f\sha\in F_0(U)$ in such a way that $f\sha$ agrees with $\varphi_2$ on
$\partial L\times[0,1]$, with $f\circ(\id\times\psi^{-1})$ on $\partial L\times[1,2]$,
and with $f$ on the complement of $Q=\partial L\times[0,2]$ in $U$. 
Then $f\mapsto f\sha$ defines a natural transformation $\nu_1$ from $F_1|\sO\tm(L;Q)$ to $F_0|\sO\tm(L;Q)$.
A natural transformation $\nu_0$ in the opposite direction, from $F_0|\sO\tm(L;Q)$ to $F_1|\sO\tm(L;Q)$, can be
defined by almost the same formula. Of course, $\varphi_2$ must be replaced by a map $\varphi_3$ from $\partial L\times[0,1]$ to
$\RR^k$ such that $\varphi_3(x,t)=\varphi_1(x,t)$ for $t\in[0,\varepsilon]$ and
such that the map from $\partial L\times[0,2]$ to $\RR^k$ defined by $(x,t)\mapsto \varphi_3(x,t)$
for $t\in[0,1]$ and $(x,t)\mapsto\varphi_0(x,t-1)$ for $t\in[1,2]$ is smooth. \newline
Now we need to show that the composition $\nu_1\nu_0\co F_0|\sO\tm(L;Q) \to F_0|\sO\tm(L;Q)$
is a weak homotopy equivalence.
For this purpose, fix $\delta$ so that $\varepsilon>\delta>0$.
Let $F^\delta_0$ be the subfunctor of $F_0$ such that $F^\delta_0(U)$ for $U\in\sO\tm(L)$ consists
of the maps $f\in F_0(U)$ which agree with $\varphi_0$ on $\partial L\times[0,\delta]\subset Q$. We will show that
$\nu_1\nu_0$ restricted to $F^\delta_0|\sO\tm(L;Q)$ is naturally homotopic to the inclusion
of $F^\delta_0|\sO\tm(L;Q)$ in $F_0|\sO\tm(L;Q)$. We need a
description of $\nu_1\nu_0$ restricted to $F^\delta_0|\sO\tm(L;Q)$. For $U$ in $\sO\tm(L;Q)$
and $f\in F^\delta_0(U)$, the map $\nu_1\nu_0(f)\co L\to \RR^k$ agrees
\begin{itemize}
\item on $\partial L\times[\frac{3}{2},2]$ with $f\circ(\id\times\psi^{-2})$;
\item on the complement of $\partial L\times[0,2]$ with $f$;
\item on $\partial L\times [0,\frac{3}{2}+\delta]$ with a map $\varphi_4$
which does not depend on $f$, where
\[
\begin{array}{rcl}
\varphi_4(x,t) & = & \varphi_0(x,t) \textup{ for }t\in [0,\delta], \\
\varphi_4(x,t) & = & \varphi_0(x,t-\frac{3}{2}) \textup{ for }t\in[\frac{3}{2},\frac{3}{2}+\delta].
\end{array}
\]

\end{itemize}
Now choose a diffeomorphism $\sigma\co [0,2]\to[0,2]$ which
maps $[0,\delta]$ diffeomorphically to $[0,\frac{3}{2}+\delta]$ and   
agrees with $\psi^2$ on $[\delta,2]$. We also require that $\sigma(t)=t$ for $t$ close to $0$ and
$\sigma(t)=t+\frac{3}{2}$ for $t$ close to $\delta$. Let $\bar\sigma\co L\to L$
be the diffeomorphism defined by $\bar\sigma(x,t)=(x,\sigma(t))$ for $(x,t)\in \partial L\times[0,2]\cong Q$
and $\bar\sigma(y)=y$ for all $y\in L\smin Q$. Then the maps
\[  f\mapsto \nu_1\nu_0(f)~, \qquad f\mapsto \nu_1\nu_0(f)\circ\bar\sigma \]
on $F^\delta_0(U)$ are homotopic, as one can see by choosing an appropriate isotopy from $\sigma$
to the identity (inducing an isotopy from $\bar\sigma$ to the identity). The map
\begin{equation} \label{eqn-firsthalfofhomotopy}  f\mapsto \nu_1\nu_0(f)\circ\bar\sigma
\end{equation}
has the following property: $\nu_1\nu_0(f)\circ\bar\sigma$ agrees with $f$ outside $\partial L\times[0,\delta]\subset Q$
while on $\partial L\times[0,\delta]$ it agrees with a map $\varphi_5$ which does not at all depend on $f$. We now
choose a path $\gamma=(\gamma_s\co\partial L\times[0,\delta]\to \RR^k)_{s\in[0,1]}$ starting with $\gamma_0=\varphi_5$ and
ending with the restriction of $\varphi_0$. We choose this so that each $\gamma_s$ agrees with $\varphi_0$ near
$\partial L\times \partial[0,\delta]$. The path induces another homotopy from~(\ref{eqn-firsthalfofhomotopy})
to the inclusion $F^\delta_0(U)\to F_0(U)$. The existence of such a path is guaranteed again by the ``large enough''
property of $\admi_\bullet$. (Note that when we selected the maps $\varphi_2$ and $\varphi_3$, we used
the assumption that the admissible maps form a dense subspace in the space of all maps satisfying the relevant
boundary conditions; but in selecting the path $\gamma$ we used the assumption that paths of admissible maps
form a dense subspace in the space of paths of all maps satisfying the relevant boundary conditions.) \newline
Now, for $U\in \sO\tm(L;Q)$, let $K$ be any compact subset of $F_0(U)$. Then, as in the proof of lemma~\ref{lem-bdrystrgwk},
it is easy to deform $K$ into a subspace of the form $F_0^\delta(U)$, for some $\delta>0$.
This means that $\nu_1\nu_0$ is homotopic to
the identity on $K$. While this does not show that $\nu_1\nu_0$ is homotopic to the identity, it is enough
to show that $\nu_1\nu_0$ is a weak homotopy equivalence, as claimed.
\qed

\proof[Sketch proof of proposition~\ref{prop-bdryindep}] Given $Q\cong\partial L\times[0,2]$ as in the
proof of lemma~\ref{lem-bdrystrgwk}, we define $Q^\varepsilon$ to be the part of $Q$ corresponding to
$\partial L\times [0,2\varepsilon]$. For $U\in \sO\tm(L)$ we let
\[  \NN_U=\{n\in \NN~|~Q^{1/n}\subset U~\}~. \]
Let $B\NN_U$ be the classifying space of $\NN_U$ as an ordered set. This is clearly contractible.
>From the proof of lemma~\ref{lem-bdryindep} we have a natural map
\[   B\NN_U\times F_0(U) \lra F_1(U) \]
which is a weak homotopy equivalence for every $U$. In addition we have the natural projection
$B\NN_U\times F_0(U)\to F_0(U)$ which is a homotopy equivalence. These two constitute our chain of natural
weak homotopy equivalences. \qed

\section{Polynomial functors, analytic functors and duality} \label{sec-polydual}
\begin{prop} Let $F$ be a good contravariant functor from $\sO\tm(L)$ to chain complexes of
abelian groups. Suppose that $F(U)$ is homotopy equivalent to a chain complex of finitely generated
\emph{free} abelian groups, bounded below, for every $U$ in $\sO\tm(L)$. Then the following are equivalent:
\begin{enumerate}
\item $F$ is polynomial of degree $\le p$~;
\item The covariant functor $DF$ defined by $DF(U)=\hom(F(U),\ZZ)$ is polynomial of degree $\le p$.
\end{enumerate}
\end{prop}

\proof We recall what it means for $F$ to be \emph{good}: it means that for $U_0\le U_1$
in $\sO\tm(L)$ such that the inclusion $U_0\to U_1$ is an isotopy equivalence (relative to
$\partial L$), the induced map $F(U_1)\to F(U_0)$ induces an
isomorphism in homology. With our assumptions, where $F(U_1)$ and $F(U_0)$ are chain complexes
of free abelian groups and bounded below, this is equivalent to saying that $F(U_1)\to F(U_0)$
is a chain homotopy equivalence. \newline
We have already given one definition of \emph{polynomial of degree $\le p$}. Namely, $F$
is polynomial of degree $\le p$ if the projection
\[  F(U) \lra \holimsub{\twosub{W\subset U}{W\in \sO p}} F(W) \]
induces an isomorphism in homology. There is another definition in terms of cubical diagrams.
Because we work in $\sO\tm(L)$ rather than $\sO(L)$, it is more technical than \cite[Def.2.2]{WeissEmb}.
Given an open set $U\in \sO(L)$ and pairwise disjoint closed subsets $C_0,\dots,C_p$ in
$U\smin\partial L$,
we can make an $S$-cube $U\smin C_\bullet$ in $\sO(L)$, where $S=\{0,1,\dots,p\}$, by
\[  T \mapsto U\smin\bigcup_{t\in T}C_t \]
for $T\subset S$. (An $S$-cube is a functor, covariant or contravariant, defined on the poset of
subsets of $S$.) To define what it means for $F$ on $\sO\tm(L)$ to be polynomial, we are
mainly interested in the case where $U\in \sO\tm(L)$ and the $C_t$ are pairwise disjoint tame co-handles.
This means that each $C_t$ is a smooth codimension $q_t$ submanifold diffeomorphic to a euclidean space, and there exists
a smooth compact codimension zero submanifold $K\subset U$ such that $\partial K$ has transverse intersections with the $C_t$~,
each intersection $K\cap C_t$ is a disk of codimension $q_t$, and the inclusion
\[ (\intr(K);\intr(K)\cap C_0,\dots,\intr(K)\cap C_p)\lra (U;C_0,\dots,C_p) \]
is isotopic to a diffeomorphism (by an isotopy in $U$ which respects the indicated submanifolds). It follows that
each $U\smin C_T$ is in $\sO\tm(L)$. We say that $F$ is polynomial of degree $\le p$ if, in this situation,
the $S$-cube of chain complexes
\[  T\mapsto F(U\smin C_T) \]
is always cartesian, where $C_T=\bigcup_{t\in T}C_t$. (For the concept of a cartesian $S$-cube,
see \cite{GoodwillieCalc2}.) This definition of \emph{polynomial of degree $\le p$} is equivalent to
the earlier one by arguments given in \cite{WeissEmb}, especially in the proofs of \cite[Thm 4.1,~Thm 5.1]{WeissEmb}.
\newline
Next, let $E$ be a \emph{covariant} functor from $\sO\tm(L)$ to cochain complexes of abelian groups. Suppose
that each $E(U)$ is homotopy equivalent to a cochain complex of finitely generated free abelian groups, bounded below.
We say that $E$ is polynomial of degree $\le p$ if, for any $U$ in $\sO\tm(L)$ and pairwise disjoint tame co-handles
$C_0,C_1,\dots,C_p$ in $U$, the $S$-cube of cochain complexes
\[  T\mapsto F(U\smin C_T) \]
is cocartesian, where $S=\{0,1,\dots,p\}$ and $C_T=\bigcup_{t\in T}C_t$. (This is our main definition of \emph{polynomial
of degree $\le p$} in the covariant setting and we studiously refrain from giving another definition which might
mention $\sO p$~.) \newline
With these definitions, the proof of our proposition turns into a triviality. Namely, $F$ with the stated finite generation and boundedness properties is polynomial
of degree $\le p$ if and only if $E=DF$ is polynomial of degree $\le p$, because the functor $\hom(\textup{--},\ZZ)$
transforms cartesian $S$-cubes of chain complexes into cocartesian $S$-cubes of cochain complexes. (Of course,
we also make use of the fact that the functor $\hom(\textup{--},\ZZ)$ is involutory on chain complexes of
finitely generated abelian groups.)
\qed

\begin{cor}
Let $F$ be a good contravariant functor from $\sO\tm(L)$ to chain complexes of
abelian groups. Suppose that $F(U)$ is homotopy equivalent to a chain complex of finitely generated
\emph{free} abelian groups, bounded below, for every $U$ in $\sO\tm(L)$. Then the following are equivalent:
\begin{enumerate}
\item $F$ is homogeneous of degree $\le p$~;
\item $DF$ is homogeneous of degree $\le p$.
\end{enumerate}
\end{cor}

\proof The functor $F$ of degree $\le p$ is homogeneous of degree $p$ if and only if
$F(U)$ is contractible for every $U\in \sO m$ where $m<p$. This can be taken as a definition.
The same definition applies in the covariant case. Therefore, $F$ is homogeneous if and only if
$DF$ is. \qed

\begin{lem} \label{lem-polydim} Let $F$ be a a good contravariant functor from $\sO\tm(L)$ to based spaces, or to chain complexes of
abelian groups. Suppose that $F(U)$ is $m$-connected for every $U\in \sO\tm(L)$. Then $T_pF(U)$
is $(m-\ell p-p)$-connected for every $U\in \sO\tm(L)$.
\end{lem}

\proof We consider first the case where the values of $F$ are based spaces (where \emph{space} means compactly generated space,
and base points are nondegenerate). The proof is by induction on $p$. For the induction beginning
we take $p=0$, in which case $T_pF(U)\simeq F(Q)$ where $Q$ is an open collar on $\partial L$. Since $F(Q)$ is
$m$-connected, this proves our claim for $p=0$. For the induction step,
we assume that $p>0$ and that $T_{p-1}F(U)$ is
$(m-\ell(p-1)-p+1)$-connected, alias $(m-\ell p-p+\ell+1)$-connected. If $m-\ell p-p+\ell+1<0$, then there is nothing
to prove for $T_pF(U)$. If $m-\ell p-p+\ell+1\ge 0$, then $T_{p-1}F(U)$ is path connected. The homotopy fiber
(over the base point) of
\[  T_pF(U) \lra T_{p-1}F(U) \]
is a homogeneous functor of degree $p$ in the variable $U$.
There is a formula \cite{WeissEmb} which describes it as a space of sections, subject to boundary
conditions, of a fibration with a
distinguished ``zero'' section on
\[   \binom{U\smin\partial L}{p}~. \]
The fiber over a configuration $S\subset U\smin\partial L$ of $p$ points is homotopy equivalent to the total
homotopy fiber of the $S$-cube
\[  R\mapsto F(N(R)) \]
where $R$ is a subset of $S$ and $N(R)$ is a tubular neighborhood of $R\cup \partial L$ in $U$.
Since $F(N(R))$ is always $m$-connected, the total homotopy
fiber of those cubes is $(m-p)$-connected. Since the configuration space has dimension $\ell p$, it follows that
the section space is $(m-\ell p-p)$-connected. Therefore the homotopy fiber of $T_pF(U)\to T_{p-1}F(U)$
is $(m-\ell p-p)$-connected, and since our connectivity estimate for $T_{p-1}F(U)$ is larger, it follows
that $T_pF(U)$ is also $(m-\ell p-p)$-connected. \newline
Now we need to look at the case where $F$ has chain complex values. This part of the proof is more
sketchy than the first. We reason that $F$ can be viewed as a functor with spectrum values.
The classification of homogeneous functors in \cite{WeissEmb} carries over to the setting of good
cofunctors on $\sO\tm(L)$ or $\sO(L)$ with values in spectra. (Fibrations over configuration
spaces have to be replaced by fibered spectra over configuration spaces.) Therefore, the proof just given
for the case of an $F$ with space values carries over to the case of an $F$ with spectrum values. \newline
We finish by explaining how chain complexes should be viewed as spectra. In the introduction, it was suggested
that the Kan-Dold equivalence can be used where necessary to view chain complexes as simplicial abelian groups, hence
via geometric realization as spaces. That suggestion needs to be made more precise.
The Kan-Dold construction $\varepsilon$ takes a chain complex $C_*$ graded
over $\ZZ$ to the simplicial abelian group $\varepsilon(C_*)$ whose set of $n$-simplices
is the set of chain maps from the cellular chain complex of $\Delta^n$ to $C_*$. It is an equivalence of
categories \emph{only} when restricted to chain complexes which are zero in negative degrees. It respects
homotopy limits in the sense that
\begin{equation} \label{eqn-KDlim}
\varepsilon\big(\!\holimsub{\alpha} C_*(\alpha) \big) \simeq \holimsub{\alpha} \varepsilon(C_*\alpha)
\end{equation}
for a functor $\alpha\mapsto C_*(\alpha)$. That formula has a weakness, though, in that $\varepsilon$ does not fully respect
properties such as \emph{not $m$-connected}. For example, it may happen that each $C_*(\alpha)$
is $3$-connected and that $\holim_\alpha C_*(\alpha)$ is $(-5)$-connected but not $(-4)$-connected. In such a case
we lose essential information by applying $\varepsilon$. The cure is to use a variant $\underline{\varepsilon}$
of the Kan-Dold construction which takes chain complexes to CW-\emph{spectra} rather than spaces. This is easy
to supply. The formula of~(\ref{eqn-KDlim}) remains valid for $\underline{\varepsilon}$.
Compared with $\varepsilon$, the construction $\underline{\varepsilon}$ has the decisive advantage that it respects
properties such as \emph{$m$-connected} and \emph{not $m$-connected}.
\qed

\medskip
We continue with a technical variant of analyticity with built-in
convergence estimates. This applies to functors with space values and to functors with
chain complex values.

\begin{defn} \label{defn-analyticrc} \cite[4.1.10]{GKW} {\rm Let $F$ be a good cofunctor on $\sO\tm(L)$. We say that
$F$ is \emph{$\rho$-analytic with excess $c$}
(where $\rho,c\in\ZZ$) if it has the following property. For $U\in\sO\tm(L)$ and $j\ge 0$ and pairwise disjoint
tame co-handles $C_0,\dots,C_j$ in $U$, of codimension $q_t<\rho$ respectively, the cube $F(U\smin C_\bullet)$
is $(c+\sum_t(\rho-q_t))$-cartesian.
}
\end{defn}

The corresponding definition in \cite{GKW} contains a very unfortunate typo or error
(where it has $k>0$ instead of $k\ge 0$, corresponding to $j\ge 0$ above). Apart from that it is
just slightly more general because tameness assumptions are
absent. As in \cite[4.2.1]{GKW} we have a convergence theorem.

\begin{thm}[Convergence] \label{thm-convergence} Suppose that $F$ is $\rho$-analytic with excess $c$, and
$U\in \sO\tm(L)$ has a tame handle decomposition relative to a collar on $\partial L$,
with handles of index $\le q$ only, where $q<\rho$.
Then the canonical map
\[  F(U) \lra T_{j-1}F(U)  \]
is $(c+j(\rho-q))$-connected, for $j>1$. Therefore the Taylor tower of $F$, evaluated at $U$,
converges to $F(U)$.
\end{thm}

(The \emph{tame handle decomposition} for $U$ is an ordinary handle decomposition for $K\subset U$, where $K\in\kappa(U)$
and the inclusion $\intr(K)\to U$ is isotopic to a diffeomorphism, relative to $\partial L$.)
\end{appendices}

\end{document}